%% file: dpg.tex
\documentclass[onefignum,onetabnum]{siamart190516}

\input{ex_shared}
\input{defs}

\input{source/trimFig}

\ifpdf
\hypersetup{
  pdftitle={Adaptive DPG for Grad-Shafranov},
  pdfauthor={Z. Peng, Q. Tang, and X.-Z. Tang}
}
\fi

\begin{document}

\maketitle

\begin{abstract}
In this work, we propose and develop an arbitrary-order adaptive
discontinuous Petrov-Galerkin (DPG) method for the nonlinear
Grad-Shafranov equation. An ultraweak formulation of the DPG scheme
for the equation is given based on a minimal residual method.  The DPG
scheme has the advantage of providing more accurate gradients compared
to conventional finite element methods, which is desired for numerical
solutions to the Grad-Shafranov equation.  The numerical scheme is
augmented with an adaptive mesh refinement approach, and a criterion
based on the residual norm in the minimal residual method is developed
to achieve dynamic refinement.  Nonlinear solvers for the resulting
system are explored and a Picard iteration with Anderson acceleration
is found to be efficient to solve the system.
Finally, the proposed algorithm is implemented in parallel on MFEM
using a domain-decomposition approach, and our implementation is
general, supporting arbitrary order of accuracy and general meshes.
Numerical results are presented to demonstrate the efficiency and
accuracy of the proposed algorithm.

\end{abstract}



\input source/sec_introduction

\input source/sec_scheme

\input source/nonlinear

\input source/sec_adaptive

\input source/sec_implementation

\input source/sec_num_0930

\section{Conclusions}
\label{sec:conclusions}

This work focuses on designing and developing an arbitrary-order
adaptive DPG method for the nonlinear Grad-Shafranov equation.  The
main focus of the work is to investigate advantages of the DPG scheme
when solving the derivatives as the magnetic field in the
Grad-Shafranov equation.  The ultraweak formulation of the DPG scheme
is developed based on a minimal residual method and the efficient
nonlinear solvers are proposed and studied.  The algorithm is
augmented with an AMR strategy to improve the efficiency of the
algorithm.  The proposed algorithm is implemented in parallel under
the framework of MFEM.

A series of numerical results is presented to verify the accuracy and
efficiency of the algorithm.  In particular, the algorithm is found to
produce more accurate numerical derivatives than some of the commonly
used finite element schemes such as CG and HDG schemes.  It is also
demonstrated that the scheme produces optimal convergence for both
the solution and derivatives for $k$-th order polynomial spaces.  The
nonlinear solver based on Anderson iteration is found to be efficient
to solve problems involving strong nonlinearity.  Finally, numerical
examples also quantitatively demonstrate the improvement of efficiency
through the AMR strategy compared with the uniform refinements.  In
conclusion, the numerical results confirm that the adaptive DPG scheme
is a good candidate for nonlinear problems when accurate derivatives
are desired, such as the Grad-Shafranov equation considered in this
work and many other problems in computational plasma physics.

\bibliographystyle{siamplain}
\bibliography{bib_peng}
\end{document}

%% file: ex_shared.tex
\usepackage{lipsum}
\usepackage{amsfonts}
\usepackage{graphicx}
\usepackage{epstopdf}
\usepackage{algorithmic}
\usepackage{tikz}
\usepackage{multirow}
\usepackage{subfigure}
\usepackage{xcolor}
\ifpdf
  \DeclareGraphicsExtensions{.eps,.pdf,.png,.jpg}
\else
  \DeclareGraphicsExtensions{.eps}
\fi

\usepackage{calc}

\usepackage{amsmath}
\DeclareMathOperator*{\argmin}{arg\,min}

\newcommand{\bogus}[1]{{}}

\newsiamremark{remark}{Remark}
\newsiamremark{rem}{Remark}
\newsiamremark{hypothesis}{Hypothesis}
\crefname{hypothesis}{Hypothesis}{Hypotheses}
\newsiamthm{claim}{Claim}

\headers{Adaptive DPG for Grad-Shafranov}{Z. Peng, Q. Tang, and X.-Z. Tang}

\title{An adaptive discontinuous Petrov-Galerkin method for the Grad-Shafranov equation\thanks{Submitted to the editors DATE.
\funding{This work was supported by the U.S. Department of
Energy through the Fusion Theory Program of the Office of Fusion
Energy Sciences, and the Tokamak Disruption Simulation (TDS)
SciDAC partnerships between the Office of Fusion Energy Science
and the Office of Advanced Scientific Computing.}}}

\author{Zhichao Peng\thanks{Department of Mathematical Sciences, Rensselaer Polytechnic Institute, Troy, NY 12180, USA.
  (\email{pengz2@rpi.edu}).}
\and Qi Tang\thanks{Los Alamos National Laboratory, Los Alamos, NM 87545, USA.
  (\email{qtang@lanl.gov}).}
\and Xian-Zhu Tang\thanks{Los Alamos National Laboratory, Los Alamos, NM 87545, USA.
  (\email{xtang@lanl.gov}).}}

\usepackage{amsopn}

\newcommand{\beq}{\begin{equation}}
\newcommand{\eeq}{\end{equation}}
\newcommand{\beqa}{\begin{eqnarray}}
\newcommand{\eeqa}{\end{eqnarray}}
\newcommand{\bit}{\begin{itemize}}
\newcommand{\eit}{\end{itemize}}
\newcommand{\bedef}{\begin{defn}}
\newcommand{\edefn}{\end{defn}}
\newcommand{\bpro}{\begin{prop}}
\newcommand{\epro}{\end{prop}}

\newcommand{\pzc}[1]{\textcolor{black}{#1}}

\newcommand{\eps}{\varepsilon}

\newcommand{\lgl}{\langle}
\newcommand{\rgl}{\rangle}

\newcommand{\wnabla}{\widetilde{\nabla}}
\newcommand{\hqn}{\widehat{q_n} }
\newcommand{\hpsi}{\widehat{\psi} }
\newcommand{\hqnh}{\widehat{q}_{n,h}  }
\newcommand{\hpsih}{\widehat{\psi}_h }
\newcommand{\bq}{   {\bf{Q} }  } 
\newcommand{\bpsi}{   {\bf{\Psi} }  } 
\newcommand{\bhqn}{   {\bf{\widehat{Q_n} } }  } 
\newcommand{\bhpsi}{   {\bf{\widehat{\Psi} } }  }

\newcommand{\bN}{   {\bf{N} }  } 
\newcommand{\bFL}{   {\bf{F_L} }  } 
\newcommand{\bL}{   {\bf{L} }  } 
\newcommand{\bu}{   {\bf{U}}      }


%% file: defs.tex
\renewcommand{\url}[1]{}
\newcommand{\citeCount}[1]{}
\newenvironment{myIndent}%
 {\list{}{\leftmargin=0.1in\rightmargin=0.1in}\item[]}%
  {\endlist}

{\noindent\textbf{Procedure~}{#1}\begin{myIndent}\em}
{\end{myIndent}}

\newcommand{\ev}{\mathbf{ e}}

\newcommand{\jv}{\mathbf{ j}}

\newcommand{\nv}{\mathbf{ n}}

\newcommand{\qv}{\mathbf{ q}}
\newcommand{\rv}{\mathbf{ r}}

\newcommand{\uv}{\mathbf{ u}}
\newcommand{\vv}{\mathbf{ v}}

\newcommand{\Bv}{\mathbf{ B}}

\newcommand{\Gv}{\mathbf{ G}}

\newcommand{\Jv}{\mathbf{ J}}

\newcommand{\Lv}{\mathbf{ L}}

\newcommand{\Qv}{\mathbf{ Q}}

\newcommand{\Uv}{\mathbf{ U}}
\newcommand{\Vv}{\mathbf{ V}}

\newcommand{\Gc}{{\mathcal G}}

\newcommand{\Jc}{{\mathcal J}}

\newcommand{\phiv}{\boldsymbol{\phi}}

\newcommand{\grad}{\nabla}

\clearpage

\newlength{\ycbTop}
\newlength{\ycbMid}%



%% file: source/trimFig.tex
\newlength{\tfwidth}
\newlength{\tfheight}
\newlength{\tfxa}
\newlength{\tfxb}
\newlength{\tfya}
\newlength{\tfyb}
\newcommand{\trimFigWithBox}[6]{%
\setlength\fboxsep{0pt}%
\setlength\fboxrule{1.0pt}
\fbox{\includegraphics[width=#2, clip, trim=#3 #4 #5 #6]{#1}}%
}
\newcommand{\trimFigNoBox}[6]{%
\setlength\fboxsep{1pt}
\setlength\fboxrule{0.0pt}
\fbox{\includegraphics[width=#2, clip, trim=#3 #4 #5 #6]{#1}}%
}
\newcommand{\trimFigHeightWithBox}[6]{%
\setlength\fboxsep{0pt}%
\setlength\fboxrule{1.0pt}
\fbox{\includegraphics[height=#2, clip, trim=#3 #4 #5 #6]{#1}}%
}
\newcommand{\trimFigHeightNoBox}[6]{%
\setlength\fboxsep{1pt}
\setlength\fboxrule{0.0pt}
\fbox{\includegraphics[height=#2, clip, trim=#3 #4 #5 #6]{#1}}%
}

\newsavebox\figBox

\newcommand{\trimw}[6]{%
\sbox\figBox{\includegraphics{#1}}
\setlength{\tfwidth}{\the\wd\figBox}
\setlength{\tfheight}{\the\ht\figBox}
\setlength{\tfxa}{\tfwidth*\real{#3}}%
\setlength{\tfxb}{\tfwidth*\real{#4}}%
\setlength{\tfya}{\tfheight*\real{#5}}%
\setlength{\tfyb}{\tfheight*\real{#6}}%
\trimFigNoBox{#1}{#2}{\tfxa}{\tfya}{\tfxb}{\tfyb}%
}

\newcommand{\trimwb}[6]{%

\sbox\figBox{\includegraphics{#1}}
\setlength{\tfwidth}{\the\wd\figBox}
\setlength{\tfheight}{\the\ht\figBox}
\setlength{\tfxa}{\tfwidth*\real{#3}}%
\setlength{\tfxb}{\tfwidth*\real{#4}}%
\setlength{\tfya}{\tfheight*\real{#5}}%
\setlength{\tfyb}{\tfheight*\real{#6}}%
\trimFigWithBox{#1}{#2}{\tfxa}{\tfya}{\tfxb}{\tfyb}%
}

\newcommand{\trimh}[6]{%
\sbox\figBox{\includegraphics{#1}}
\setlength{\tfwidth}{\the\wd\figBox}
\setlength{\tfheight}{\the\ht\figBox}
\setlength{\tfxa}{\tfwidth*\real{#3}}%
\setlength{\tfxb}{\tfwidth*\real{#4}}%
\setlength{\tfya}{\tfheight*\real{#5}}%
\setlength{\tfyb}{\tfheight*\real{#6}}%
\trimFigHeightNoBox{#1}{#2}{\tfxa}{\tfya}{\tfxb}{\tfyb}%
}

\newcommand{\trimhb}[6]{%

\sbox\figBox{\includegraphics{#1}}
\setlength{\tfwidth}{\the\wd\figBox}
\setlength{\tfheight}{\the\ht\figBox}
\setlength{\tfxa}{\tfwidth*\real{#3}}%
\setlength{\tfxb}{\tfwidth*\real{#4}}%
\setlength{\tfya}{\tfheight*\real{#5}}%
\setlength{\tfyb}{\tfheight*\real{#6}}%
\trimFigHeightWithBox{#1}{#2}{\tfxa}{\tfya}{\tfxb}{\tfyb}%
}

%% file: source/sec_introduction.tex
\section{Introduction}

The magnetohydrodynamic (MHD) equilibrium is critical in many
applications of plasma systems, such as magnetic confinement fusion.
Steady-state fusion reactor operation requires that a MHD equilibrium
is reached and sustained.  In computational plasma physics, an
efficient and accurate solver for the MHD equilibrium serves as the
basis for linear stablity analysis and nonlinear simulations of plasma
transport and off-normal events such as tokamak disruptions.  For
example, the confining magnetic fields computed from a MHD equilibrium
solver can be used as a starting point for the relativistic
Fokker-Planck equation in the runaway electron
study~\cite{guo2017models}, or a starting point to evolve a nonlinear
MHD simulation~\cite{jorek}.  For more details on the importance to have
accurate and efficient equilibrium solvers, see the discussions
in~\cite{lee2015ecom} for instance.

The axisymmetric MHD equilibrium is governed by the Grad-Shafranov
equation, which is an elliptic equation with a nonlinear source
term~\cite{grad1958hydromagnetic,shafranov1958magnetohydrodynamical}.  There are two types of equilibrium
problems, the so-called fixed boundary and free boundary equilibria.
In the fixed boundary equilibrium problem, it is assumed that the
separatrix of the MHD equilibrium is known and a Dirichlet boundary
condition is imposed at the separatrix for the Grad-Shafranov
equation. (In practice, a closed flux surface just inside the
separatrix, the so-called $q_{95}$ surface, is used for the
computational boundary.) The free boundary problem considers the case
of an unbounded domain when separatrix is not known and is in fact to
be determined self-consistently from both the plasma current and
electrical currents in the external field coils. Thus, a free boundary
solver needs extra capabilities, such as \pzc{a boundary} integral term
that addresses the effect of a far-field boundary and numerical
routines to locate the separatrix in a given trial solution.  As an
example, see~\cite{heumann2017finite} for a finite-element-based free
boundary solver.  In this work, we focus on the fixed boundary problem
and develop a fast and high-order finite-element-based solver on
general meshes.  The fixed boundary solver is a necessary building
block for a free boundary solver.

In recent years, the fixed boundary solvers have drawn significant
attentions due to its importance in magnetic confinement fusion.  Many
numerical schemes have been explored and developed. These include
hybridizable discontinuous Galerkin
methods~\cite{sanchez2019hybridizable, sanchez2019adaptive}, spectral
elements~\cite{howell2014solving, palha2016mimetic}, boundary integral
approaches~\cite{pataki2013fast, lee2015ecom}, Hermite finite
element~\cite{huysmans1991isoparametric, lutjens1996chease} and so on.
Some of them, such as~\cite{sanchez2019adaptive}, also considered
adaptive mesh refinement (AMR) approaches to minimize numerical errors
and improve the efficiency of the algorithm.  The derivative of the
solutions determines physically important quantities such as magnetic
fields for computing the charged particle
trajectory. Hence, other than efficiency and accuracy of the
algorithms, many of the previous works also focus on minimizing errors
in the derivatives of the solutions, including using bicubic elements
or introducing auxiliary variables for its derivatives. Constructing
accurate derivatives and effective AMR will be the two primary
objectives of the current work.

Among many available high-order schemes, we are interested in using
discontinuous Petrov Galerkin (DPG) methods for the Grad-Shafranov
equation.  The DPG method proposed in~\cite{demkowicz2010class,
  demkowicz2011class} enjoys the following properties that are
particularly attractive for our purpose. First, with an ultra-weak
formulation and suitable choice of the discrete test space, the DPG
method can provide high order accurate approximation to both solutions
and the first derivatives~\cite{gopalakrishnan2014analysis}. Second,
the DPG method provides a natural built-in error estimator for
adaptivity based on a numerical residual
~\cite{demkowicz2012class,carstensen2014posteriori}, while a standard
AMR method relies on the calculation of local numerical fluxes or
physical features like high order derivatives.

Next we briefly review the DPG methodology. More details on the
formulations can be found in Section
\ref{sec:dpg_minimal_residual}. The DPG method is a Petrov-Galerkin
finite element method, which uses on-the-fly computed optimal test
functions ~\cite{demkowicz2010class,
  demkowicz2011class}. Discontinuous or broken test space is used so
that these optimal test functions can be locally computed on each
element. Besides the perspective of optimal functions, the DPG method
can be equivalently seen as a minimal residual method
\cite{demkowicz2014overview,demkowicz2017discontinuous,keith2017discrete,chan2014dpg,roberts2015discontinuous}. The
DPG method preserves a discrete version of the inf-sup condition and
provides the best error approximation for an energy norm defined in
the trial space \cite{demkowicz2014overview}.
Another attractive feature of the method is that it provides a natural
built-in error
estimator~\cite{demkowicz2012class,carstensen2014posteriori}. AMR and
$hp$-adaptivity
algorithms~\pzc{\cite{demkowicz2012class,chan2014dpg,roberts2015discontinuous,petrides2017adaptive}}
are designed based on this estimator. DPG methods for nonlinear
problems have also been considered. For instance, in
\cite{chan2010class,chan2014dpg,roberts2015discontinuous}, the authors
first linearize nonlinear problems and then apply the DPG method to
the linearized problems. \pzc{In~\cite{nagaraj20193d}, a Picard iteration is applied to solve DPG method for a nonlinear optical problem. } In~\cite{moro2012hybridized}, from the
perspective of minimal residual method, a hybridized DPG method is
designed and directly applied to nonlinear fluid problems, and a
Newton method is used for the nonlinear solver. In~\cite{bui2014pde},
a partial-differential-equation-constrained optimization method is
considered as a nonlinear solver for the DPG method. Recently,
\pzc{a posteriori error analysis} is also generalized to the nonlinear problem
in~\cite{carstensen2018nonlinear}.

In this work, we would like to take full advantage of the DPG
\pzc{methodology} as described above and develop an efficient and accurate
solver for the Grad-Shafranov equation, which produces not only
high-order accurate solution but also high-order accurate derivatives.
The major contributions of this work include: (1) we extended the DPG
scheme to the nonlinear Grad-Shafranov equation using the formulation
based on a minimal residual method; (2) nonlinear solvers to the
DPG formulation are considered, including both Newton's methods and a
Picard iteration method, of which the Picard iteration method
accelarated with Anderson acceleration and algebraic multigrid
preconditioners is found to be more efficient; (3) an AMR approach is
developed based on the residual norm coming along with the DPG scheme;
(4) the algorithm is implemented in parallel on MFEM ({\tt mfem.org}),
and the solver supports arbitrary order \pzc{accuracy} and general meshes for
complex geometry.

The remainder of the paper is structured as follows.  The MHD
equilibrium and Grad-Shafranov equation are briefly described in
Section~\ref{sec:governing}.  The formulation of the DPG scheme is
given in Section~\ref{sec:scheme}.  Our discussion starts
with a review of the abstract DPG formulation based on minimal
residual method.  The ultraweak formulation is then described
 for the Grad-Shafranov equation.  The details on
discretizations such as its matrix-vector form follow.  In
Section~\ref{sec:nonlinear}, we focus on efficient nonlinear solvers
for the resulting matrix-vector form of the DPG scheme.
In Section~\ref{sec:adaptive}, we discuss adaptive mesh refinement,
followed by details of the implementation in
Section~\ref{sec:implementation}.  Finally, numerical results are
presented in Section~\ref{sec:num}.

\section{Governing equation}
\label{sec:governing}
In this section, we briefly introduce the MHD equilibrium and
Grad-Shafranov equation, and then define the fixed-boundary problem
that is solved by the DPG scheme in this work.

A MHD equilibrium for magnetically confined plasma means the force balancing 
\[
    \jv\times\Bv = \grad p,
\]
where $\Bv$ is the magnetic field, $\jv$ is the plasma current density
that satisfies Ampere's law
\[
    \mu_0 \jv = \grad \times \Bv,
\]
with $\mu_0$ the magnetic permeability, and $p$ the plasma pressure.
If further assuming the problem is axisymmetric (e.g., the equilibria
in tokamaks) and considering the equilibrium in $(r,z)$-coordinate,
the Grad-Shafranov equation for the magnetic flux function $\psi$ can
be written as
\begin{align*}
    r \frac \partial {\partial r} \left (\frac 1 r \frac{\partial
      \psi}{\partial r}\right) + \frac{\partial^2 \psi}{\partial z^2}
    = -\mu_0 r^2 \frac{dp}{d\psi} - I \frac{dI}{d\psi},
\end{align*}
where $I$ is a function of $\psi$ and it is associated with the
toroidal part of the magnetic field.  The magnetic field satisfies
\[
    \Bv = \frac 1 r \grad \psi \times \hat e_\theta + \frac{I} r \hat
    e_\theta
\]
Note that the quantities of interest in practice are the derivatives
of $\psi$ (corresponding to the magnetic field) and its second
derivatives (corresponding to the current).  For more details on the
Grad-Shafranov equation and plasma equilibrium, the readers are
referred to plasma textbooks such as~\cite{jardin2010computational}.

Define a gradient operator $\wnabla = (\partial_r, \partial_z)^T$.  In
this work, we consider the Grad-Shafranov equation with the fixed
boundary condition and the problem is rewritten as
\begin{alignat}{3}
\label{eqn:gs}
&\wnabla\cdot \left({\frac{1}{r}\wnabla \psi }\right) = - \frac{F(r,z,\psi)}{r}, && \qquad x\in\Omega,\notag\\
&\psi = \psi_D, && \qquad x\in \partial \Omega,
\end{alignat}
where $\Omega\subset \mathbb{R}^2$ is the physical domain with a
Lipschitz boundary $\partial\Omega$.  In tokamaks, the physical domain
$\Omega$ corresponds to the cross section of the device.  The source
term $F(r, z, \psi)$ depends on $p$ and $I$, and in practice both of
them are functions of $\psi$, which are from either experimental
measurements or theoretical design.  Therefore, the problem is
nonlinear through the source term and constructing efficient nonlinear
solvers is one of the main focuses of the current work.

%% file: source/sec_scheme.tex
\section{Formulation of the DPG method}
\label{sec:scheme}
We begin with a quick review of the formulation of the DPG method as a
minimal residual method for a general abstract nonlinear problem.
Based on this abstract formulation, we will show some details of an
ultraweak DPG formulation for the Grad-Shafranov equation.

\subsection{Abstract DPG method as a minimal residual method}
\label{sec:dpg_minimal_residual}
Consider an abstract weak formulation for a general nonlinear problem.
Let $U$ denote the trial space and $V$ the test space, both of
which are Hilbert spaces.  The weak formulation for a nonlinear
problem is to find $u\in U$ such that,
\begin{align}
b_N(u,v) = l(v), \qquad \forall \, v\in V,
\end{align}
where $l(\cdot): V\rightarrow \mathbb{R}$ is a linear form and
$b_N(u,v): U\times V\rightarrow \mathbb{R}$ is nonlinear for $u\in U$
while linear for $v\in V$.  A residual operator can be defined as
\begin{align}
r(u,v) = b_N(u,v) - l(v).
\end{align}
Let $V'$ denote the dual space of $V$.  Since $b_N(u,v)$ is linear for
$v$, 
\pzc{ based on the Riesz representation theory, there exists a linear operator $\mathcal{B}: U\rightarrow V'$ such that
 $\forall u\in U$,  $<\mathcal{B}u,v>_{V'\times V} = b_N(u,v)$, $\forall v\in V$.}  It is also convenient to define the inner products
$(\cdot,\cdot)_U$ and $(\cdot,\cdot)_V$ as the inner products in $U$
and $V$, respectively.

As pointed out in~\cite{moro2012hybridized, bui2014pde, keith2017discrete,carstensen2018nonlinear} for nonlinear problems, the DPG method is
equivalent to a minimal residual method.  In this work, we focus on
this \pzc{interpretation} of the DPG method in both the discussion and
implementation.  
\pzc{Here, we briefly review the minimal residual interpretation, and more discussions can be found in \cite{moro2012hybridized,demkowicz2014overview,chan2014dpg,roberts2014dpg,demkowicz2017discontinuous,keith2017discrete}.}
Suppose $U_h\subset U$ and $V_h\subset V$ to be the finite dimensional
discrete trial and test spaces.  The DPG method can be defined as
looking for $u_h\in U_h$ such that
\begin{align}
\label{eq:min_res}
u_h &= \argmin_{w_h\in U_h} || \pzc{\mathcal{B}{\omega_h}  }- l ||_{V'}^2
\notag\\
&= \argmin_{w_h\in U_h} \left( \sup_{v_h\in V_h,\; v_h\neq 0} \frac{|b_N(w_h,v_h) - l(v_h)|^2 }{||v_h||^2_V}\right)
\notag\\
&= \argmin_{w_h\in U_h} \left( \sup_{v_h\in V_h,\; v_h\neq 0} \frac{|r(w_h,v_h)|^2 }{||v_h||^2_V}\right).
\end{align}
Let $\{e_{v_i} \}_{i=1}^N$ be a basis of the discrete test space $V_h$
and ${\bf{e_v} }:= ( e_{v_1},\dots, e_{v_N} )^T$, then for $\forall
v_h \in V$, it can be expanded as $v_h = {\bf{v} }^T {\bf{e_v} }=
\sum_{i=1}^N v_i \, e_{v_i}$. Since both $l(v)$ and $b_N(u,v)$ are
linear with respect to the argument $v$, we have
\begin{align*}
    r(w_h,v_h) = r\left(w_h, \sum_{i=1}^N v_i \, e_{v_i} \right) = \sum_{i=1}^N v_i r\left(w_h,e_{v_i} \right) = {\bf{v} }^T {{\rv(w_h)} },
\end{align*}
where ${{\rv(w_h)} } := \left(r(w_h,e_{v_1}),\dots,r(w_h,e_{v_N})\right)^T$. 
Then, we note that
\begin{align}
\label{eq:matrix_form_dual_norm}
 \frac{|r(w_h,v_h)|^2 }{||v_h||^2_V} = \frac{|{\bf{v}}^T{{\rv(w_h} })|^2}{ {\bf{v} }^T G {\bf{v} } },
\end{align}
where the matrix $G\in \mathbb{R}^{N\times N}$ is defined as
\begin{align}
\label{eqn:def_G}
G_{ij} := \left(e_{v_i}, e_{v_j} \right)_V.
\end{align}
$G$ is often called a Gram matrix, and here it can be interpreted as a
mass matrix corresponding to the inner product $(\cdot,\cdot)_V$.
 

{The Gateaux derivative of ~\eqref{eq:matrix_form_dual_norm} with
  respect to ${\bf{v}}$ is}
\begin{align}
    \frac{2}{ ({\bf{v} }^TG{\bf{v} })^2 }\Big[ ({\bf{v}}^T\rv(w_h))({\bf{v} }^TG{\bf{v} })\rv(w_h)-({\bf{v}}^T\rv(w_h))^2G{\bf{v} }\Big].
\end{align}
It is easy to see when $\vv$ satisfies
\[
{\bf{v} } = \frac{ {\bf{v} }^TG{\bf{v} }   }{{\bf{v} }^T {\rv(w_h) } } G^{-1}{\rv(w_h)},
\]
the Gateaux derivative is $0$, and the supremum of
$\frac{|r(w_h,v_h)|^2 }{||v_h||^2_V}$ is obtained as
 \begin{align}
  \sup_{v_h\in V_h,\; v_h\neq 0} \frac{|r(w_h,v_h)|^2 }{||v_h||^2_V}
  = {{\rv(w_h)} }^T G^{-1} {{\rv(w_h)} }.
 \end{align}
Therefore, the DPG method \eqref{eq:min_res} is equivalent to find 
\begin{align}
\label{eq:min_res2}
u_h = \argmin_{w_h\in U_h} \Big( {{\rv(w_h)} }^T G^{-1} {{\rv(w_h)} } \Big).
\end{align}
Taking the Gateaux derivative of \eqref{eq:min_res2} gives
\begin{align}
\label{eqn:dpg}
J^T(u_h) G^{-1} {\bf{r(u_h)} } = 0, 
\end{align}
where $J(u_h)$ is the Jacobian matrix of  the residual ${{\rv(u_h)} }$ with repect to $u_h$.
Finally, the DPG method based on the minimal residual method is given by~\eqref{eqn:dpg},
which will be the focus of the current work. \pzc{More details for the minimal residual interpretation of the DPG methodology can be found in \cite{moro2012hybridized,keith2017discrete}.}

The original DPG scheme was defined by computing optimal test functions on the fly, 
see~\cite{demkowicz2011class} for instance. 
The DPG scheme in~\cite{demkowicz2011class} is equivalent to the form~\eqref{eqn:dpg}
derived from the minimal residual method. 
We refer the readers to~\cite{moro2012hybridized} for details of its equivalence proof. 
In summary, the methodology of the DPG method uses special ``optimal'' test functions, which \pzc{preserves} a discrete inf-sup condition with \pzc{a discrete inf-sup constant of the same order as the original continuous inf-sup constant} and \pzc{guarantees} the ``optimal convergence'' results,
see \cite{demkowicz2011class,demkowicz2014overview} for instance.


Finally, consider the special case of  linear problems. 
Let  $\{e_{u_i} \}_{i=1}^M$ be the basis of the discrete trial space $U_h$ and ${\bf{e_u} }:=(e_{u_1},\dots,e_{u_M})^T$. For $\forall u_h\in U_h$, it can be expanded as $u_h = {\bf{u} }^T {\bf{e_u} }$. For a linear problem when $b_N(u,v)$ becomes a bilinear form, 
one can verify that~\eqref{eqn:dpg} becomes 
\[
\pzc{B^T G^{-1} B {\bf{u} } = B^T G^{-1}{\bf{\ell} }, }
\] 
where the matrix $B\in \mathbb{R}^{N\times M}$ defined by $B_{ij} =
b(e_{u_j}, e_{v_i})$ and ${\bf{\ell} }$ is a constant vector
determined by operator $l(\cdot)$. Hence, for linear problems, the DPG
method always results in a symmetric positive definite linear system,
while nonlinear problems do not have such a property, which will be
discussed in Section~\ref{sec:nonlinear}.

\subsection{An ultraweak DPG formulation for the Grad-Shafranov equation}
We present the details of the DPG formulation for the
problem~\eqref{eqn:gs}.  The motivation to use the DPG scheme is to
obtain a more accurate approximation to the gradient of $\psi$, and we
therefore consider the ultraweak formulation to
discretize~\eqref{eqn:gs}.  Define an auxiliary variable of the vector
\[
\mathbf{q}:= -\frac{\wnabla \psi} r,
\]
and rewrite \eqref{eqn:gs} into its first order form,
\begin{subequations}
\label{eqn:gs_1st}
\begin{alignat}{3}
&r{\bf{q} } +\wnabla \psi = 0,\quad  && \qquad x\in\Omega,\\
&\wnabla \cdot {\bf{q} } = \frac{1}{r}F(r,z,\psi), &&\qquad x\in\Omega,\\
& \psi = \psi_D,  && \qquad x\in \partial\Omega.
\end{alignat}
\end{subequations}
Let $\Omega_h=\{K\}$ be a partition of the physical domain $\Omega$,
where $K\in \Omega_h$ are disjoint elements. Let $\partial \Omega_h =
\{ \partial K, K\in T_h\}$ denote its skeleton (edge) and $\Gamma_h =
\partial \Omega_h \cap \partial \Omega $ be the set of the edges on
the physical boundary. Define $(\cdot,\cdot)_K$ and
$\lgl\cdot,\cdot\rgl_{\partial K}$ as the standard $L^2$ inner product
of $L^2(K)$ and $L^2(\partial K)$, respectively.  Let
$(\cdot,\cdot)_{\Omega_h} = \sum_{K\in\Omega_h }(\cdot,\cdot)_K$,
$\lgl\cdot,\cdot\rgl_{\partial\Omega_h} = \sum_{\partial
  K\in\partial\Omega_h }\lgl \cdot,\cdot\rgl_{\partial K}$ and
$\lgl\cdot,\cdot\rgl_{\Gamma_h} = \sum_{\partial K\in \Gamma_h }\lgl
\cdot,\cdot\rgl_{\partial K}$.

We generalize the DPG method for the linear Poisson problem
in~\cite{gopalakrishnan2014analysis} to the nonlinear Grad-Shafranov
equation.  
\pzc{Define the trace spaces 
\[
    H^{-\frac{1}{2} }(\partial \Omega_h) :=\{\hqn: \, \exists \, {\bf{q} }\in H(\textrm{div};\Omega)\;\;\text{such that}\;\; \hqn = {\bf{q} }\cdot n|{\partial_K},\;\;\forall K\in\Omega_h\},
\]
\[
    H^{\frac{1}{2}}(\partial \Omega_h) := \{\hpsi: \, \exists\, \psi\in H^1(\Omega)\;\;\textrm{such that}\;\; \hpsi = \psi|{\partial_K},\;\;\forall K\in\Omega_h\},
\]
and the broken Sobolev spaces
\[
H(\textrm{div},\Omega_h):=\{\phiv\in L^2({\Omega_h}): \,\,\phiv|K\in H(\textrm{div},K),\;\forall K\in \Omega_h\},
\]
\[
H^{1}(\Omega_h):=\{\tau\in L^2({\Omega_h}): \,\,\tau|K\in H^1(K),\;\forall K\in \Omega_h\}.
\]
}
The trial space $U$ and the test space $V$ are \pzc{chosen} as
\begin{align*}
&U = \left(L^2(\Omega) \right)^2\times L^2(\Omega) \times H^{-\frac{1}{2}}(\partial \Omega_h) \times H^{\frac{1}{2} }(\partial \Omega_h),\\
&\pzc{V =H(\textrm{div},\Omega_h) \times H^{1}(\Omega_h)}.
\end{align*}
The ultraweak form associated with~\eqref{eqn:gs_1st} is to find ${\bf{u} } = ({\bf{q} },\psi,\hqn, \hpsi)\in U$ such that for $\forall \, {\bf{v} }= (\phiv,\tau) \in V$
\begin{subequations}
\label{eqn:dpg_gs}
\begin{align}
\label{eqn:dpg_q}
&(r{\bf{q} },\phiv)_{\Omega_h} -\left( \psi, \wnabla\cdot \phiv \right)_{\Omega_h} + \lgl \hpsi, \nv \cdot \phiv \rgl_{\partial {\Omega_h} }=0,\\
\label{eqn:dpg_divq}
-&({\bf{q} },\wnabla \tau)_{\Omega_h} + \lgl \hqn, \tau\rgl_{\Omega_h}= \left(\frac{F(r,z,\psi)} r, \tau\right)_{\Omega_h},\\
\label{eqn:dpg_bc}
& \lgl \hpsi, \tau \rgl_{\Gamma_h} = \lgl \psi_D, \tau\rgl_{\Gamma_h}.
\end{align}
\end{subequations} 
The source term $F(r,z,\psi)$ can be rewritten as the summation of a nonlinear part and a linear part $F(r,z,\psi)=F_N(r,z,\psi)+F_L(r,z)$. Based on the weak form \eqref{eqn:dpg_gs}, we define 
\begin{subequations}
\label{eqn:dpg_operator}
\begin{align}
\label{eqn:dpg_l}	   
l({\bf{v} }) &= \left(\frac{F_L(r,z)} r, \tau\right)_{\Omega_h},
\\
\label{eqn:dpg_b}
b_N({\bf{u} },{\bf{v} })& = (r{\bf{q} },\phiv)_{\Omega_h} -\left( \psi, \wnabla\cdot \phiv\right)_{\Omega_h} + \lgl \hpsi, \nv\cdot \phiv \rgl_{\partial {\Omega_h} } 
\notag\\
	  & -({\bf{q} },\wnabla \tau)_{\Omega_h} + \lgl \hqn, \tau\rgl_{\Omega_h}- \left(\frac{F_N(r,z,\psi)} r, \tau\right)_{\Omega_h}.
\end{align}
\end{subequations}
In addition, we define the test norm $||\cdot||_V$ as
\begin{align}
\label{eqn:dpg_test_norm}
||{\bf{v} }||_V^2 := || (\phiv,\tau)||_V^2 = ||\phiv||^2 + || \wnabla\cdot \phiv||^2 + ||\tau||^2 + ||\wnabla \tau||^2,
\end{align}
and let $(\cdot,\cdot)_V$ be its corresponding inner product.

Finally, we need to determine the discrete trial space $U_h$ and the
discrete test space $V_h$. Let $P^k(K)$ and $P^k(\partial K)$ be the
space of polynomials with degree at most $k$ on the element $K$ and
its edge $\partial K$. 
The discrete trial space $U_h$ and the discrete test space $V_h$ are respectively chosen as
\begin{alignat*}{2}
    U_h^k &= \Big\{ &&u_h = (\qv_h,\psi_h, \hqnh, \hpsi_h): \, \qv_h|K \in  \left(P^k(K)\right)^2,\; \psi_h|K \in P^k(K),\; \\
		&  &&\hqnh|\partial K \in P^k(\partial K)\cap H^{-\frac{1}{2} }(\partial \pzc{\Omega_h}),\; \hpsih|\partial K\in P^{\pzc{k+1} }(\partial K)\cap H^{\frac{1}{2} }(\partial \pzc{\Omega_h})
	  \Big\},
\\	  
V_h^{k,s} &= \Big\{ && v_h= (\phiv_h,\tau_h): \, \phiv_h|K \in \left(P^{k+s}(K) \right)^2,\; \tau_h|K\in P^{k+s}(K)\Big\}, \quad s\geq 2.	  
\end{alignat*}
According to \cite{carstensen2014posteriori}, the key to prove the
convergence of the DPG scheme is to define a Fortin-type operator
$\Pi: V\rightarrow V_h=V_h^{k,s}$ such that $b_N(w_h, (\mathbf{I}-\Pi)v)=0$
for $\forall w_h\in U_h, v\in V$. As proved in
\cite{gopalakrishnan2014analysis}, $s\geq 2$ is a sufficient condition
for the existence of such a Fortin-type operator $\Pi$, which further
guarantees the convergence of the scheme.  We therefore choose $s\ge2$
in this work.  \pzc{A parameter study for the Fortin operator is carried out in \cite{nagaraj2017construction}.} Finally, based on the given $U_h$, $V_h$,
\eqref{eqn:dpg_operator} and \eqref{eqn:dpg_test_norm}, we are ready
to determine $G$, $J(u_h)$ and ${{\rv(u_h) } }$ in $\eqref{eqn:dpg}$
and define our DPG scheme.

\pzc{
\begin{rem}
The choice of the test norm $||\cdot||_V$ is not unique. An attractive alternative is the adjoint graph norm \cite{zitelli2011class,demkowicz2014overview,petrides2019adaptive}:
\begin{align}
||{\bf{v} }||_{V,\textrm{adjoint} }^2 := || (\phiv,\tau)||_{V,\textrm{adjoint}}^2 = ||r\phiv-\wnabla \tau||^2 + || \wnabla\cdot \phiv||^2+ ||\phiv||^2 + ||\tau||^2 ,
\end{align}
which delivers a quasi-optimal error estimate in $L^2$ sense. In \cite{demkowicz2011analysis,gopalakrishnan2014analysis}, our current choice $||\cdot||_V$ is also proved to provide a quasi-optimal error estimate for linear elliptical problems.
\end{rem}
}

\subsection{Matrix-vector form of the DPG scheme}
In this section, the details of the DPG scheme are presented as a matrix-vector form. 
The following notations will be used throughout this section. 
Let 
\begin{align*}
    U_h \supset \{e_{u_j} \}_{j=1}^M := \{ {\ev }^{\qv}_{j_{\qv}}\}^{M_{\qv}}_{j_{\qv} =1} \times \{ e^{\psi}_{{j_\psi}}\}^{M_\psi}_{j_\psi =1}\times \{ e^{\hqn}_{j_{\hqn} } \}^{M_{\hqn} }_{j_{\hqn} =1}\times  \{ e^{\hpsi}_{j_{\hpsi}}\}^{M_{\hpsi} }_{j_{\hpsi} =1}
\end{align*}
and
\begin{align*}
    V_h \supset \{e_{v_i}\}_{i=1}^N :=  \{  {{\ev} }^{\phiv}_{i_{\phiv}}\}^{N_{\phiv}}_{i_{\phiv} =1}\times \{ e^{\tau}_{{i_\tau}}\}^{N_\tau}_{i_\tau =1} 
\end{align*}
be the bases of the discrete trial space $U_h$ and the test space $V_h$,  respectively.
Note that ${\bf{e} }^{\qv}_j$ and ${\bf{e} }^{\phiv}_i$ are vector-valued, while other basis functions are scalar-valued.
Define $\ev_\uv^{\xi}:=(e^\xi_1,\dots, e^\xi_{M_\xi} )^T$ and $\ev_\vv^{\eta} := (e^\eta_1,\dots, e^\eta_{N_\eta} )^T$, where $\xi =\qv, \, \psi, \, \hqn, \, \hpsi$ and $\eta = \phiv,\,\eta$. Then, $q_h$, $\psi_h$, $\hqnh$ and $\hpsih$ can be rewritten as
\begin{align*}
\qv_h = \bq^T \ev_\uv^\qv, \qquad \psi_h = {\bf{\Psi} }^T \ev_\uv^\psi,\qquad \hqnh = \bhqn^T \ev_\uv^{\hqn},\qquad \hpsih=\bhpsi^T\ev_\uv^{\hpsi}.
\end{align*}
Define $\bu:= (\bq^T, \bpsi^T, \bhqn^T, \bhpsi^T)^T$,
and on each element $\Qv$ contains two components corresponding to the vector $\qv$.
Now, we are ready to rewrite the ultraweak formulation~\eqref{eqn:dpg_gs} into its matrix-vector form and obtain the residual as
 \begin{align}
\label{eqn:residual}
{{\rv(u_h)} } &= b_N(u_h,v_h) - l(v_h) \notag \\
&= B_L \Uv - B_N(\Uv) -\bFL \notag\\
& = \left(\begin{matrix}
	M_r   & \wnabla_h & 0             & T_{\hpsi} \\
    \textrm{div}_h &  0           & T_{\hqn}  & 0
	 \end{matrix}  \right)
	 \left(
	 \begin{matrix}
	 \bq \smallskip\\
	 \bpsi \smallskip\\
	 \bhqn \smallskip\\
	 \bhpsi  
	 \end{matrix}
	 \right)
-\left(\begin{matrix}
0 \\
\bN(\bpsi)
\end{matrix}\right)
-\left(\begin{matrix}
0 \\
\bL
\end{matrix}\right),
\end{align}
where the block matrices $M_r$, $\wnabla_h$, $T_{\hpsi}$, $\textrm{div}_h$ and $T_{\hqn}$ are defined as 
\begin{alignat*}{4}
&(M_r)_{ij} := \left( {\bf{e} }^{\qv}_j,{\bf{e} }^{\phiv}_i\right)_{\Omega_h},\quad 
&& (\wnabla_h)_{ij}:= -\left(e^{\psi}_j, \wnabla \cdot {\bf{e} }^{\phiv}_i\right)_{\Omega_h}, & \\ 
& (T_{\hpsi})_{ij} := \lgl  e^{\hpsi}_j,\nv\cdot {\bf{e} }^{\phiv}_i\rgl_{\partial \Omega_h}, \quad
&& (\textrm{div}_h)_{ij} := -\left( {\bf{e} }^{\qv}_j, \wnabla e^\tau_i\right)_{\Omega_h},\quad
 (T_{\hqn})_{ij} := \lgl  e^{\hqn}_j,e^\tau_i\rgl_{\partial \Omega_h}, 
\end{alignat*} 
and the vectors $\bN(\bpsi)$ and $\Lv$ are defined as
\begin{align*}
 (\bN(\bpsi) )_i := \left(F_N(r,z,\psi_h)/r, e_i^\tau\right)_{\Omega_h},\qquad 
 \bL _i := \left(F_L(r,z)/r, e_i^\tau\right)_{\Omega_h}. &
\end{align*}
The Jacobian matrix of ${{\rv(u_h)} }$ can be therefore derived as
\begin{align}
\label{eqn:jacobi}
J(u_h) 
&=
\left(\begin{matrix}
	M_r   & \wnabla_h & 0             & T_{\hpsi} \\
    \textrm{div}_h &  -D_N(\psi_h) & T_{\hqn}  & 0
	 \end{matrix}  \right)
\end{align}
where the block matrix $D_N(\psi_h)$ is defined as $(D_N(\psi_h) )_{ij} :=  \left( \frac{\partial \left(F_N(r,z,\psi_h)/r\right)}{\partial \psi} e_j^\psi, e_i^\tau\right)_{\Omega_h}$ and $\frac{\partial \left(F_N(r,z,\psi_h)/r\right)}{\partial \psi}$ is the derivative of $F_N(r,z,\psi_h)/r$ with respect to its third argument.
Recall the definitions of the Gram matrix $G$ in \eqref{eqn:def_G} and the test norm $||\cdot||_V$ in \eqref{eqn:dpg_test_norm}, and we get 
\begin{align}
\label{eqn:gram}
G=\left(\begin{matrix}
G_V & 0 \\
0    & G_S
\end{matrix}\right),
\end{align}
where the block matrices are defined as 
\begin{align*}
    (G_V)_{ij}&=\left( {\bf{e} }^{\phiv}_j,{\bf{e} }^{\phiv}_i\right)_{\Omega_h}+\left( \wnabla\cdot{\bf{e} }^{\phiv}_j,\wnabla\cdot{\bf{e} }^{\phiv}_i\right)_{\Omega_h},\\
(G_S)_{ij}&=\left( e^\tau_j,e^\tau_i\right)_{\Omega_h}+\left( \wnabla e^\tau_j,\wnabla e^\tau_i\right)_{\Omega_h}.
\end{align*}
Substituting \eqref{eqn:residual}, \eqref{eqn:jacobi} and \eqref{eqn:gram} into \eqref{eqn:dpg}, we further rewrite the matrix-vector form of the DPG scheme as
\begin{align}
\label{eqn:dpg-mat-vec}
J^T(\bu) G^{-1} \Big[B_L \Uv - B_N(\bu) -\bFL \Big]=0.
\end{align}
The above system is the final formulation of the DPG scheme for the Grad-Shafranov equation.
The  system is nonlinear due to the source term and thus it requires special care to solve it efficiently. In Section~\ref{sec:nonlinear} we will focus on exploring its nonlinear solvers.

%% file: source/nonlinear.tex
\section{Nonlinear solvers}
\label{sec:nonlinear}
In this work, we have considered two types of nonlinear solvers, Newton's method and Picard iteration equipped with Anderson acceleration, for the matrix-vector form ~\eqref{eqn:dpg-mat-vec}.
Some advantages and disadvantages are found during the numerical experiments,
which will be addressed in this section.
All the nonlinear solvers discussed in this work 
are implemented through PETSc SNES solvers~\cite{balay2019petsc}.

\subsection{Newton's method}
The most popular method for nonlinear problems is Newton's method,
which is the first type of solvers we have investigated during this
work.  The challenge to implement a Newton's method for the
system~\eqref{eqn:dpg-mat-vec} is to assemble its Jacobian matrix. It
is easy to see that the Jacobian of the system~\eqref{eqn:dpg-mat-vec}
is given by
\begin{align}
    \Jc = J^T(\bu) G^{-1} J(\Uv) + 
\left(\begin{matrix}
	0   & 0 & 0   & 0 \\
    0 &  \mathcal{H} \rv_2 & 0  & 0 \\
	0   & 0 & 0   & 0 \\
	0   & 0 & 0   & 0 
     \end{matrix}  \right).
     \label{eqn:exact-j}
\end{align}
Here $\rv_2$ stands for the second component in the residual $\rv(u_h)$ and
$\mathcal{H}$ is the Hessian defined by $- \frac{\partial D_N}{\partial \bpsi}$ where $D_N$
is the matrix associated with the nonlinear source term.
To implement the Jacobian matrix $\Jc$, one has to implement a Hessian-vector product associated
with $\mathcal{H}$, which is not trivial especially in finite element packages such as MFEM.
Instead, we avoid  assembling the Hessian matrix by using the Jacobian-free Newton-Krylov (JFNK) 
method that estimate $\Jc$ through a finite difference approximation of 
\[
    \Jc\approx\frac{\Gv(\Uv+\epsilon\Vv)-\Gv(\Uv)}\epsilon,
\]
where the nonlinear function in this case is $\Gv(\Uv) = J^T(\bu) G^{-1} \Big[B_L \Uv - B_N(\bu) -\bFL \Big]$.
Without a proper preconditioner,  the JFNK solver is typically very inefficient. 
Through looking at its exact form of~\eqref{eqn:exact-j}, it seems natural to provide a preconditioner of $J^T(\bu) G^{-1} J(\Uv)$ for JFNK 
which only ignore one diagonal term related to $\mathcal{H}$.
A block Jacobi preconditioner is further provided to accelerate the inversion of $J^T(\bu) G^{-1} J(\Uv)$ (the block Jacobi preconditioner is identical to the preconditioner for $\Jv(\Uv)G^{-1}B_L$
that will be discussed in the next section and we therefore skip the discussion here).
A GMRES solver for the approximated $\Jc$ is further used in the outer iteration to guarantee the convergence of the full nonlinear solver.

Such an approach works fine for problems that have very weak nonlinearity,
while it fails to converge for hard problems when the nonlinear source term is strong.
The ignorance of the block matrix $\mathcal{H} \rv_2$ appears to be critical in those problems.
For instance, for the nonlinear problems presented in Section~\ref{sec:num}, 
the JFNK solver with the above preconditioning strategy fails,  while the same solver works well
for the linear problems presented there.
Therefore, we further explore another approach that is based on Picard iteration,
which is discussed in the next section.

\subsection{Anderson acceleration method}
The second type of nonlinear solvers we have investigated is the Picard iteration with the Anderson
acceleration (Anderson mixing) method~\cite{anderson1965iterative,walker2011anderson}.
We present the algorithm of Anderson acceleration that is used in this work
\begin{algorithm}[H]
\caption{Anderson acceleration  solving nonlinear fixed point problem $\bu-A(\bu)=0.$ }
\label{alg:buildtree}
\begin{algorithmic}
\STATE{Given initial guess $\bu^0$. Set the initial residual ${\bf{R} }_0$ as ${\bf{R} }_0 =\bu_0 - A(\bu_0)$.}
\STATE{Set $\bu_1 = (1-\lambda_0) U_0 +\lambda_0 A(\bu_0)$, ${\bf{R} }_1 =\bu_1 - A(\bu_1)$  and $k=1$.}
\WHILE{$\Big(||r_k||_2\geq\max(rtol ||r_0||_2,atol)$ and $||x_k - x_{k-1}||\geq stol ||x_{k-1}||\Big)$}
\STATE{Set $m_k=\min\{k,m\}$.}
\STATE{Find $(\alpha_0^{k},\dots,\alpha_{m_k}^k)$ to be the solution of the constrained minimization problem:
$$(\alpha_0^{k},\dots,\alpha_{m_k}^k)= \arg \min \left\| \sum_{i=0}^{m_k}\alpha_i^k {\bf{R} }_{k-i}\right \|_2^2,\; \text{such that}\; \sum_{i=0}^{m_k}\alpha_{i}^k=1.$$
}
\STATE{Set $\bu_{k+1} = (1-\lambda_k)\sum_{i=0}^{m_k} \alpha^k_i \bu_{k-i} + \lambda_k \sum_{i=0}^{m_k} \alpha^k_i A({\bu}_{k-i})$.}
\STATE{Calculate the residual: ${\bf{R} }_{k+1} = \bu_{k+1} - A(\bu_{k+1} )$ and $k:=k+1$.	}
\ENDWHILE
\end{algorithmic}
\end{algorithm}
Here $\lambda_k$ is determined by a cubic backtracking line search algorithm \cite{dennis1996numerical}. 

Note that the system~\eqref{eqn:dpg-mat-vec} is not in the form of a fixed point problem,
and hence some transformations are needed before applying the Anderson acceleration.
We have investigated two types of rewritten systems.
A straightforward rewritten form of~\eqref{eqn:dpg-mat-vec}  compatible with the Anderson acceleration leads to the following system 
$$\Uv - \left\{\Uv-J^T(\bu) G^{-1} \Big[B_L \Uv - B_N(\bu) -\bFL \Big]\right\}=0.$$ 
However, we found that the nonlinear solver for this system converges  poorly even with the Anderson acceleration. 
Instead, the following nonlinear fixed-point problem is used in our implementation
\begin{align}
\label{eqn:dpg-num}
\bu - \left(J^T(\bu) G^{-1} B_L\right)^{-1}\Big[ J^T(\bu) G^{-1}(B_N(\bu)-\bFL ) \Big]=0.
\end{align}
It is straightforward to see that ~\eqref{eqn:dpg-num} is equivalent with~\eqref{eqn:dpg-mat-vec}.
Note that a similar strategy was applied in \cite{sanchez2019hybridizable,sanchez2019adaptive}. 
According to \cite{brune2015composing}, this reformulation strategy can be viewed as a ``preconditioning'' procedure for the nonlinear solver.

At each iteration, we need to invert the block matrix of $J^T(\Uv)G^{-1}B_L$. 
Using \eqref{eqn:residual}, \eqref{eqn:jacobi} and \eqref{eqn:gram}, with some work, we found that
\begin{align*}
&J^T(\bu) G^{-1}B_L \\
=&\left(\begin{matrix}
	M_r   & \wnabla_h & 0             & T_{\hpsi} \\
    \textrm{div}_h &  -D_N(\psi_h) & T_{\hqn}  & 0
	 \end{matrix}  \right)^T
\left(\begin{matrix}
G_V & 0 \\
0    & G_S
\end{matrix}\right)^{-1}
 \left(\begin{matrix}
	M_r   & \wnabla_h & 0             & T_{\hpsi} \\
    \textrm{div}_h &  0 & T_{\hqn}  & 0
     \end{matrix}  \right)
\\
=&
\left(\begin{matrix}
M_r^T G_V^{-1} M_r + \textrm{div}_h ^TG_S^{-1} \textrm{div}_h                   &   M_r^T G_V^{-1}\wnabla_h             & \textrm{div}_h^T G_S^{-1} T_{\hqn}           & M_r^T G_V^{-1} T_{\hpsi}\\
\wnabla^T_h G_V^{-1} M_r - D_N(\psi_h)^T G_S^{-1} \textrm{div}_h   & \wnabla_h^TG_V^{-1}\wnabla_h      &-D_N(\psi_h)^T G_S^{-1}T_{\hqn} & \wnabla^T G_V^{-1}T_{\hpsi}\\ 
T_{\hqn}^T G_S^{-1} \textrm{div}_h                                                   & 0                                                       & T_{\hqn}^T G_S^{-1} T_{\hqn }     & 0\\
T_{\hpsi}^T G_V^{-1} M_r                                                     & T_{\hpsi}^T G_V^{-1}\wnabla_h       & 0                                                   & T_{\hpsi}^T G_V^{-1} T_{\hpsi}
\end{matrix}\right)
\end{align*}
Note that the block matrix is not symmetric due to the appearance of the term $-D_N(\psi_h)$,
and the flexible GMRES solver \cite{saad1993flexible} is applied to invert this block matrix.
Since the final form of $J^T(\Uv)G^{-1}B_L$ is rather complicated, it is necessary to provide a few words on some details of matrix assembling when implementing the scheme.
In our implementation, we assemble most of small blocks in the beginning of the simulation
except for two blocks associated with $-D_N(\psi_h)$.
During each iteration of the nonlinear solver, we only need to update two small block matrices 
resulting from the nonlinear source term.
The inversion of the Gram matrices, $G_V$ and $G_S$,  is implemented exactly by inverting the corresponding small
matrix on each element during assembling.
This is sufficient since our test space is discontinuous.
The rest of matrix assembling uses standard approaches.

The iterative linear solver to invert $J^T(\Uv)G^{-1}B_L$ is still inefficient without preconditioning.
Since this inversion needs to be performed at each iteration of the nonlinear solver, 
it is required that such an iterative linear solver must be efficient.
To improve the efficiency of this linear solver, we provide a block Jacobian preconditioner associated with $J^T(\Uv)G^{-1}B_L$
\begin{align*}
P  =
\left(\begin{matrix}
P_{11} &              &              &\\
           & P_{22}   &              &\\
           &               & P_{33}  &\\
           &               &              & P_{44}
\end{matrix}\right) 
\end{align*}
where the blocks are given by
\begin{alignat*}{3}
& P_{11}:= M_r^T G_V^{-1} M_r + \textrm{div}_h^TG_S^{-1} \textrm{div}_h, \qquad
 &&  P_{22}:=\wnabla_h^TG_V^{-1}\wnabla_h, \\
    &   P_{33}:= T_{\hqn}^T G_S^{-1} T_{\hqn }, \qquad   
    &&  P_{44}:= T_{\hpsi}^T G_V^{-1} T_{\hpsi}
\end{alignat*}
In the preconditioning stage, we inverted each block of $P_{ii}$ with algebraic multigrid preconditioners provided by HYPRE.
We found it is sufficient to use the standard algebraic multigrid method (AMG) \cite{ruge1987algebraic} with one V-cycle to invert the blocks of $P_{11}$,  $P_{22}$ and $P_{44}$. The AMS method \cite{kolev2009parallel} with one V-cycle is applied to invert
$P_{33}$. 
Note that the block $P_{33}$ is corresponding to  $\hqnh$, which is associated with the trace of a $H(\textrm{div})$ space. As an AMG method designed for $H(\textrm{div})$ and $H(\textrm{curl})$ spaces based on the idea of auxiliary space method \cite{hiptmair2007nodal}, the AMS method improves the efficiency over the standard AMG method for those spaces. 
Hence, compared with the standard AMG method, 
we found that the AMS method is a better choice to invert the block $P_{33}$.

%% file: source/sec_adaptive.tex
\section{Adaptivity}
\label{sec:adaptive}
An adaptive mesh refinement (AMR) approach is developed in our implementation to improve the efficiency of our DPG solver.
The basic idea of AMR is to refine a mesh dynamically based on some estimate of local errors
in the current solution. For the steady state problem such as the case we consider in this work,
AMR can improve the convergence of nonlinear solvers for some hard problems compared with
a uniformly refined mesh, since 
a significant amount of work is focused on the region where errors are large in an adaptive mesh thanks to the error estimator.
One advantage of the DPG method in the context of AMR is that it comes with a natural error estimator, while other numerical schemes typically rely on certain estimates, such as the  Zienkiewicz-Zhu error estimator which estimates the errors using local numerical fluxes or
error estimators through Richardson extrapolation.
For instance, in \cite{carstensen2014posteriori} a natural posteriori error estimator based on the residual $||B_u - l ||_{V'}$ is derived for 
 the DPG method solving linear problems. 
 In \cite{carstensen2018nonlinear} the posteriori error analysis is generalized to nonlinear problems and an AMR algorithm is presented for the nonlinear problem.

Here we first present the details of the error estimator for the DPG method and then describe the details of an AMR algorithm based on it.

\subsection{Error estimator}
Recall \eqref{eq:min_res2}, and the DPG method seeks $u_h$ such that
\begin{align}
  \pzc{u_h=\argmin_{w_h\in U_h} || \mathcal{B}{\omega_h} - l ||_{V'}^2.}
\end{align}
A posteriori error estimate for DPG method solving linear problems is given in \cite{carstensen2014posteriori} as
\begin{align}
\pzc{C_1 ||\mathcal{B}{u_h}-l ||_{V'} \leq || {\bf{u} } - u_h||_U \leq C_2||\mathcal{B}{u_h}-l ||_{V'}  +C_3 \textrm{osc}(l),}
\end{align}
where $C_1$, $C_2$ and $C_3$ are constants independent of the mesh. The term $\textrm{osc}(l)$ does not depend on the numerical solution $u_h$, and it is a high order term with respect to the mesh size $h$, determined by the approximation properties of the finite element spaces. With regularity assumptions, \cite{carstensen2018nonlinear} generalizes the results to nonlinear problems.
As shown in Section \ref{sec:dpg_minimal_residual}, 
\begin{align}
\pzc{
|| \mathcal{B}{u_h} - l ||^2_{V'}
}= \sup_{v_h\in V_h,\; v_h\neq 0} \frac{|r(u_h,v_h)|^2 }{||u_h||^2_V}
={\bf{r(u_h)} }^T G^{-1} {\bf{r(u_h)} }.
\end{align}
Using the Riesz representation, we can find $\eps_{u_h}\in V$ such that 
$$(\eps_{u_h},v)_V =b_N(u_h,v)-l(v)= r(u_h,v) = (\pzc{\mathcal{B}{u_h} }-l)(v),\forall v\in V. $$
Furthermore, $|| \eps_{u_h}||_{V} = ||\pzc{\mathcal{B}{u_h}}-l||_{V'}=\sqrt{{\bf{r(u_h)} }^T G^{-1} {\bf{r(u_h)} }}$.
Then, we define the local error estimator $E_K$ on each element $K$ as
\begin{align}
\label{eqn:estimator}
E_K := || \eps_{u_h} ||_{V(K)}. 
\end{align}
The corresponding total error estimator for the whole domain is defined as
\begin{align}
\label{eqn:estimator_total}
E_\textrm{total}:=\sqrt{ \sum_{K\in \Omega_h} ||\eps_{u_h}||^2_{V(K)}}.
\end{align}

\subsection{Adaptive strategy}
We only consider conforming  $h$-adaptive refinement with a fixed
polynomial order in this work, although 
the non-conforming AMR is recently developed in MFEM~\cite{vcerveny2019non}.
We use the following mesh-refinement strategy in the implementation.
After each nonlinear solver for the DPG scheme,  
an updated solution $\uv_h$ for the Grad-Shafranov equation is obtained and its residual norm on each element can be estimated using~\eqref{eqn:estimator}.
In our implementation, for a given element $K$, it will be refined if all the following three conditions hold:
\begin{equation*}
\left\{ 
\begin{alignedat}{3}
&(1) \; &&\eps_K >\textrm{atol}_\textrm{amr},  \\
&(2) \; &&\eps_K >\theta_{\textrm{max}} \, \max_{K'\in\Omega_h} \eps_{K'},  \\
&(3)\; &&\eps_K>\theta_{\textrm{total} } \, E_{\textrm{total} }/\sqrt{N_{mesh}},
\end{alignedat} \right.
\end{equation*}
where $\theta_\textrm{max}\in [0,1)$, $\theta_\textrm{total}\in [0,1)$ and $\textrm{atol}_\textrm{amr}\geq 0$ are some given thresholds.
In summary, we refine the element $K$ based on the criterion that its local error is less than a predetermined tolerance $\rm atol_\textrm{amr}$,
and its estimated error is relatively larger compared to the errors in other elements.
Finally, the stoppage criterion for the AMR iteration is when the total number of elements is larger than some
upper bound, or {no element satisfies the three conditions above}.

%% file: source/sec_implementation.tex
\section{Implementation}
\label{sec:implementation}

The described DPG algorithm for the Grad-Shafranov equation is implemented under the framework of c++ library MFEM~\cite{mfem_url}. 
It would be appropriate to say a few words about our implementation under MFEM.

The algorithm is implemented in parallel using a standard domain decomposition method.
All the vectors and small block matrices use the parallel distributed data structure provided by the package  
and its communication between sub-domains is based on the message-passing interface (MPI).
All the linear and nonlinear solvers described in this work are implemented through PETSc~\cite{balay2019petsc} and 
some of the small matrices are preconditioned with HYPRE \cite{falgout2002hypre} algebraic multigrid preconditioners to improve efficiency.
All the matrix assemblies are performed under MFEM to minimize the communication cost between MFEM and PETSc.

Our implementation of the DPG algorithm is general, 
supporting arbitrary order of accuracy and general meshes including triangular, quadrilateral and high-order
curvilinear meshes, taking full advantage of the capabilities from the package.
However, our current implementation only supports conforming AMR, but the implementation still offers refining and dynamic load-balancing,
which will be performed after the mesh is updated in each AMR iteration.  
As a result, in Section~\ref{sec:num}, we mainly focus on the numerical results on triangular meshes and
 its conforming adaptive mesh refinement to verify the scheme.
Another reason to focus on triangular meshes  is due to the limitation of mesh generations. 
Currently, for all the problems with complex geometry, we generate triangular meshes using the Gmsh software. 
To the best of our knowledge, curvilinear meshes with complex geometry obtained from a mesh generator are not fully supported by MFEM.
Nevertheless, the focus of the current work is on the DPG scheme, AMR and efficient nonlinear solvers,
while mesh generators are well beyond the scope of this work. 
Therefore, it is sufficient to use triangular meshes generated by Gmsh for the numerical tests.

\paragraph{Reproducibility} The implementation of the DPG scheme for the Grad-Shafranov equation, as well as the numerical examples presented in this work, is freely available as an MFEM fork at {\tt github.com/ZhichaoPengMath/mfem}.

%% file: source/sec_num_0930.tex
\section{Numerical tests}
\label{sec:num}
We demonstrate the performance of our DPG scheme through a set of numerical examples. 
Several examples are presented to study the accuracy of the DPG scheme, including two linear examples and one nonlinear example.
The linear examples are further used to compare the performance of the DPG scheme with other finite element methods.
We then consider a nonlinear case involving slightly more complicated geometry.
Finally, several nonlinear examples are presented to study the performance of AMR.

In practice,  we use the trace of the Raviart-Thomas space \cite{raviart1977mixed,nedelec1980mixed} for the space corresponding to the third component of $U_h^k$. For the discrete test space $V_h^{k,s}$,  $s=2$ is chosen throughout the numerical study.
Unless otherwise noted, the results presented here use the Picard iteration with Anderson acceleration as the nonlinear solver.

\subsection{Accuracy tests}
\label{sec:accuracy}
Two numerical examples are presented to demonstrate the accuracy of the DPG scheme. 
The first example considers  the linear Solov’ev profiles from \cite{pataki2013fast,cerfon2010one}. This test is further used to compare the 
 accuracy for both the solution $\psi$ and its derivative $\wnabla{\psi}/r$ solved by our DPG method, 
 conventional continuous Galerkin (CG) method and the hybridized discontinuous Galerkin (HDG) method. The second example 
 considers a manufactured solution with a nonlinear source. We use it to demonstrate the accuracy of our scheme in nonlinear problems.

\subsubsection{Linear tests}
\label{sec:accuracy_linear}
The linear Grad-Shafranov equation of 
$$\pzc{r}\wnabla\cdot \left({\frac{1}{r}\wnabla \psi }\right) = r^2$$
 has an exact solution \cite{pataki2013fast,cerfon2010one} in the form of
$$ \psi(r,z) = \frac{r^4}{8} + d_1 + d_2 r^2 + d_3(r^4- 4 r^2z^2),$$
where the parameters $d_1$, $d_2$ and $d_3$ determine the reasonable plasma cross section. 

{
\newcommand{\figWidth}{6.0cm}
\newcommand{\trimfig}[2]{\trimh{#1}{#2}{.2}{.2}{.2}{.25}}
\begin{figure}[htb]
\begin{center}
\begin{tikzpicture}[scale=1]
  \useasboundingbox (0,.5) rectangle (12.,6.);  
  \draw(-.7,-.6) node[anchor=south west,xshift=-4pt,yshift=+0pt] {\trimfig{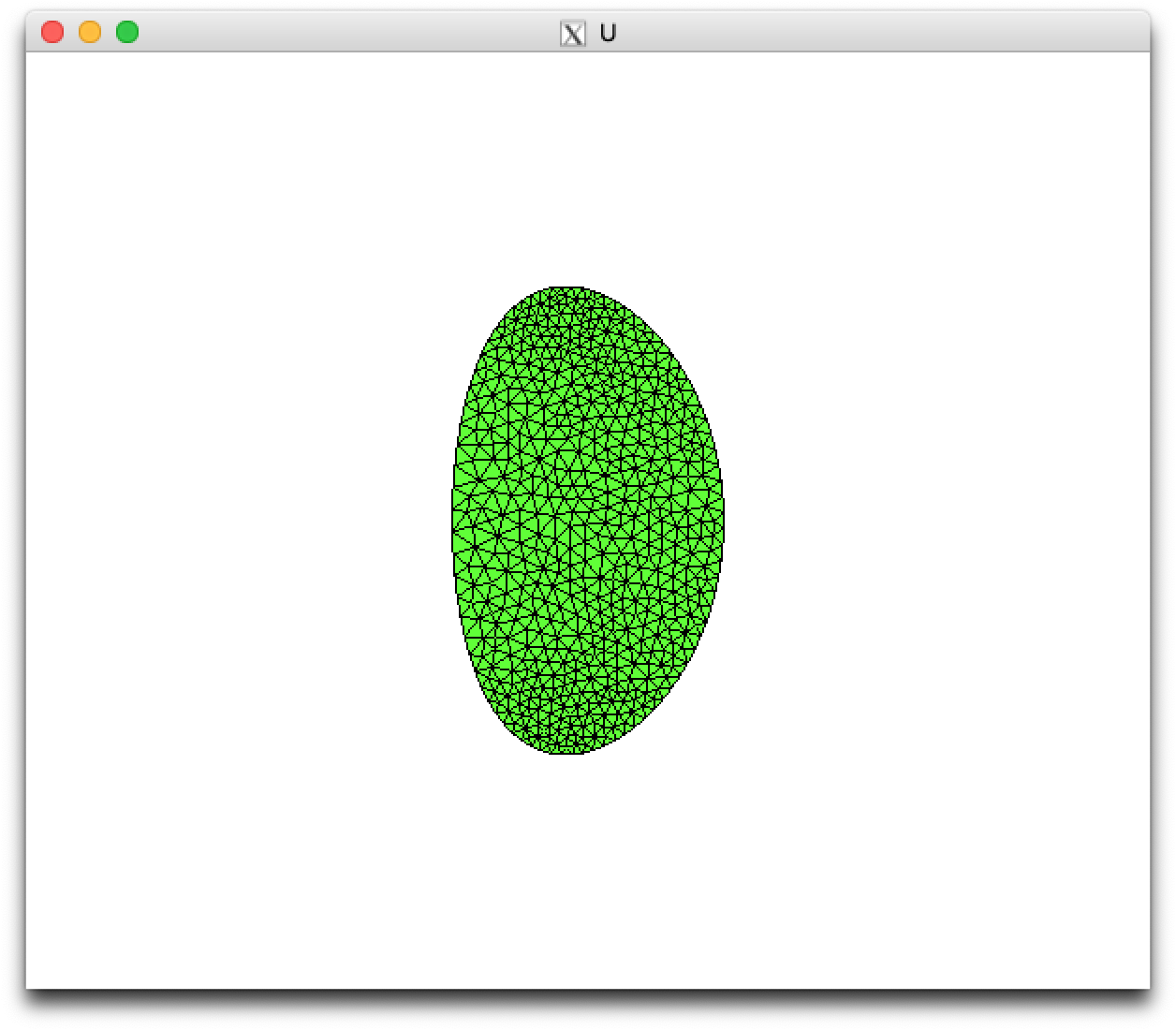}{\figWidth}};
  \draw(5, -.6) node[anchor=south west,xshift=-4pt,yshift=+0pt] {\trimfig{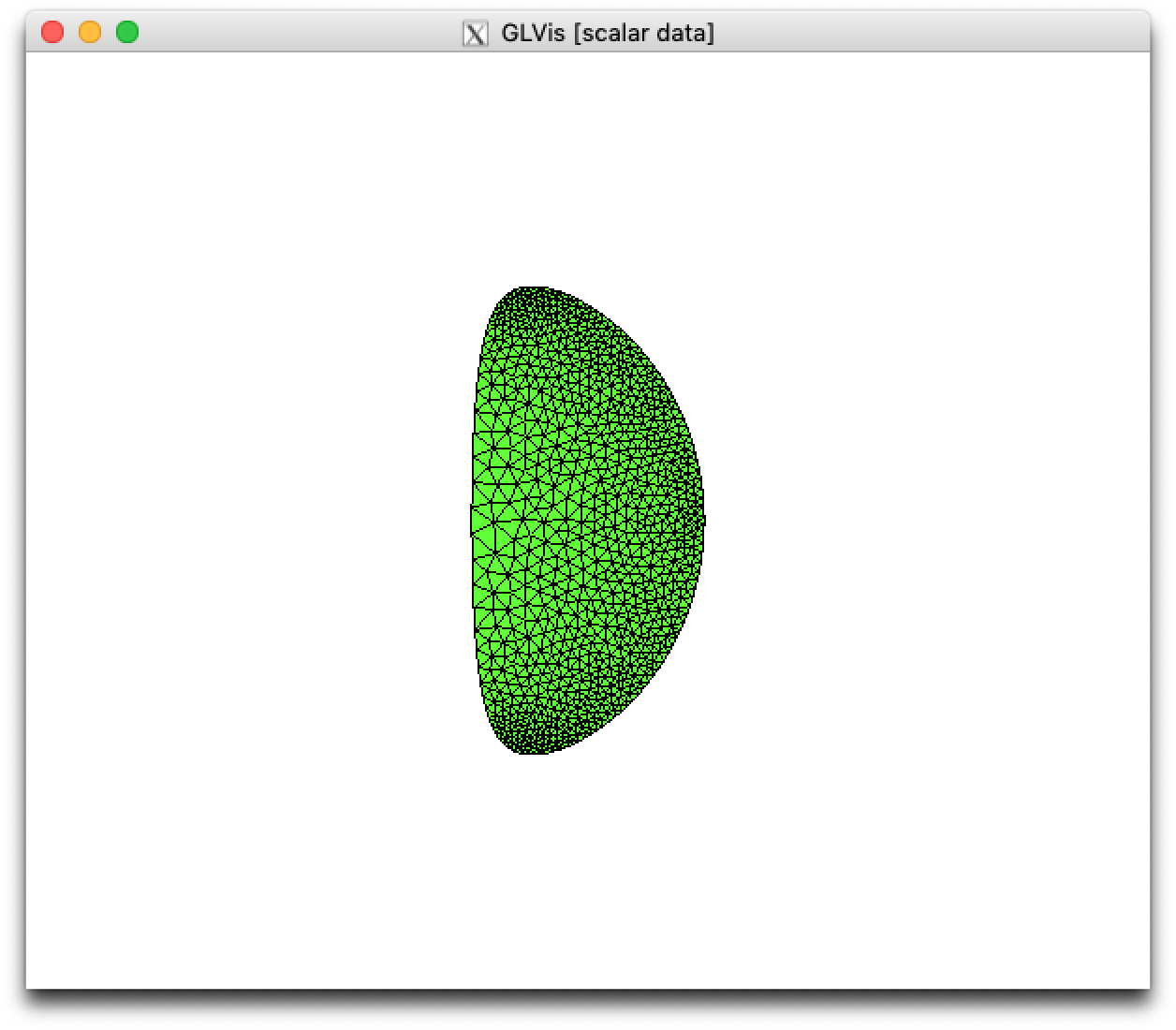}{\figWidth}};
  \draw(3.3,5.6) node[draw,fill=white] {\footnotesize $\Gc^{(1)}$ of ITER};
   \draw(9,5.6) node[draw,fill=white] {\footnotesize $\Gc^{(1)}$ of NSTX};
%
\end{tikzpicture}
\end{center}
  \caption{Computational domains and initial meshes of $\Gc^{(1)}$ for two linear accuracy tests.
  The meshes are generated through Gmsh by describing the computational boundaries.
  }
\label{fig:accuracy_mesh}
\end{figure}
}

In \cite{pataki2013fast} it has been described how to use the inverse
aspect ratio $\eps$, the elongation $\kappa$ and the triangularity
$\delta$ to determine the parameters $d_1$, $d_2$ and $d_3$.  We adopt
the same manufactured solutions in the accuracy test. The parameters
are determined by the following linear system
\begin{align*}
\left(\begin{matrix}
1 & (1+\eps)^2 & (1+\eps)^4\\
1 & (1-\eps)^2  & (1-\eps)^4 \\
1 & (1-\delta\eps)^2 & (1-\delta\eps)^4 - 4(1-\delta\eps)^2\kappa^2\eps^2
\end{matrix}\right)
\left(\begin{matrix}
d_1\\
d_2\\
d_3\\
\end{matrix}\right)
=-\frac{1}{8}\left(\begin{matrix}
(1+\eps)^4\\
(1-\eps)^4\\
(1-\delta\eps)^4\\
\end{matrix}\right).
\end{align*}
Two tests are considered in this case, including the ITER-like \cite{aymar2002iter} configuration of $\eps=0.32$, $\kappa=1.7$, $\delta=0.33$  and NSTX-like configuration \cite{sabbagh2001equilibrium} of $\eps=0.78$, $k=2$, $\delta = 0.35$. The computational domains and initial meshes are presented in Figure~\ref{fig:accuracy_mesh}.
The meshes are generated using Gmsh by describing the computational domains using the exact solutions. 
In particular, the computational boundaries are described by Chebyshev nodes along the $r$ direction and $z$ coordinates are then determined by solving $\psi = 0$.
 The numerical solutions are shown in Figure~\ref{fig:accuracy1} and the presented results include $\psi$ and its two derivatives of $\psi_r/r$ and $\psi_z/r$.
The results show a good agreement with their exact solutions.

 {%
\newcommand{\drawContour}[7]{%
\begin{scope}[#1]
\draw(0.0,-.3) node[anchor=south west,xshift=-4pt,yshift=+0pt] {\trimfig{picture/#2}{\figWidth}};
\draw(1.5,5.) node[draw,fill=white,anchor=west,xshift=2pt,yshift=0pt] {\small #3};
\draw(.5,5.) node[draw,fill=white,anchor=west,xshift=2pt,yshift=0pt] {\small #5};
\begin{scope}[xshift=2pt]
  \draw (\xcb,\ycb) node[anchor=south west,xshift=0.cm,yshift=.5cm,rotate=-90] {\trimfigcb{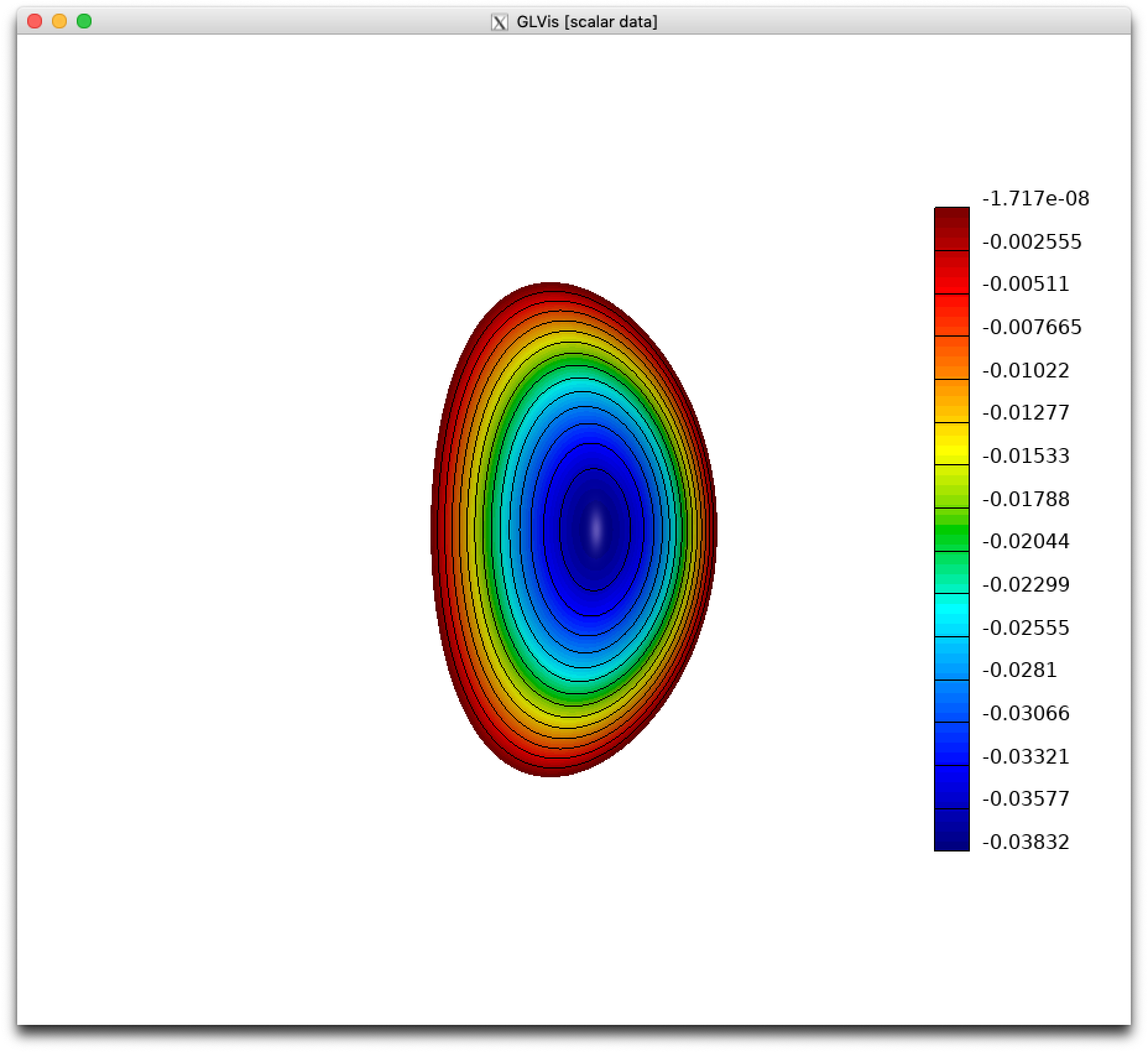}{\cbWidth}{\cbHeight}};
  \draw (.8,0) node[anchor=north,xshift=+3pt,yshift=-2pt] {\scriptsize $#6$};
  \draw (4.4,0) node[anchor=north,xshift=+0pt,yshift=-2pt] {\scriptsize $#7$};
\end{scope}
\end{scope}
}
\newcommand{\cbWidth}{.3cm}
\newcommand{\cbHeight}{4cm}
\newcommand{\xcb}{.5cm}
\newcommand{\ycb}{-.2cm}
\setlength{\ycbTop}{\ycb+\cbHeight}
\setlength{\ycbMid}{\ycb+\cbHeight*\real{.5}}
\newcommand{\trimfigcb}[3]{\includegraphics[width=#2, height=#3, clip, trim=26.5cm 4cm 5cm 4cm]{#1}}
%
\newcommand{\figWidth}{5cm}
\newcommand{\yb}{6.75}
\newcommand{\trimfig}[2]{\trimw{#1}{#2}{.25}{.25}{.2}{.2}}
\begin{figure}[htb]
\begin{center}
\resizebox{13cm}{!}{
\begin{tikzpicture}[scale=1]
  \useasboundingbox (0.0,.5) rectangle (16.,12);  
 \drawContour{xshift= 0.0cm,yshift=6.75cm}{test1_psi.png}{ITER}{}{$\psi$}{-0.0358}{0};
 \drawContour{xshift=5.25cm,yshift=6.75cm}{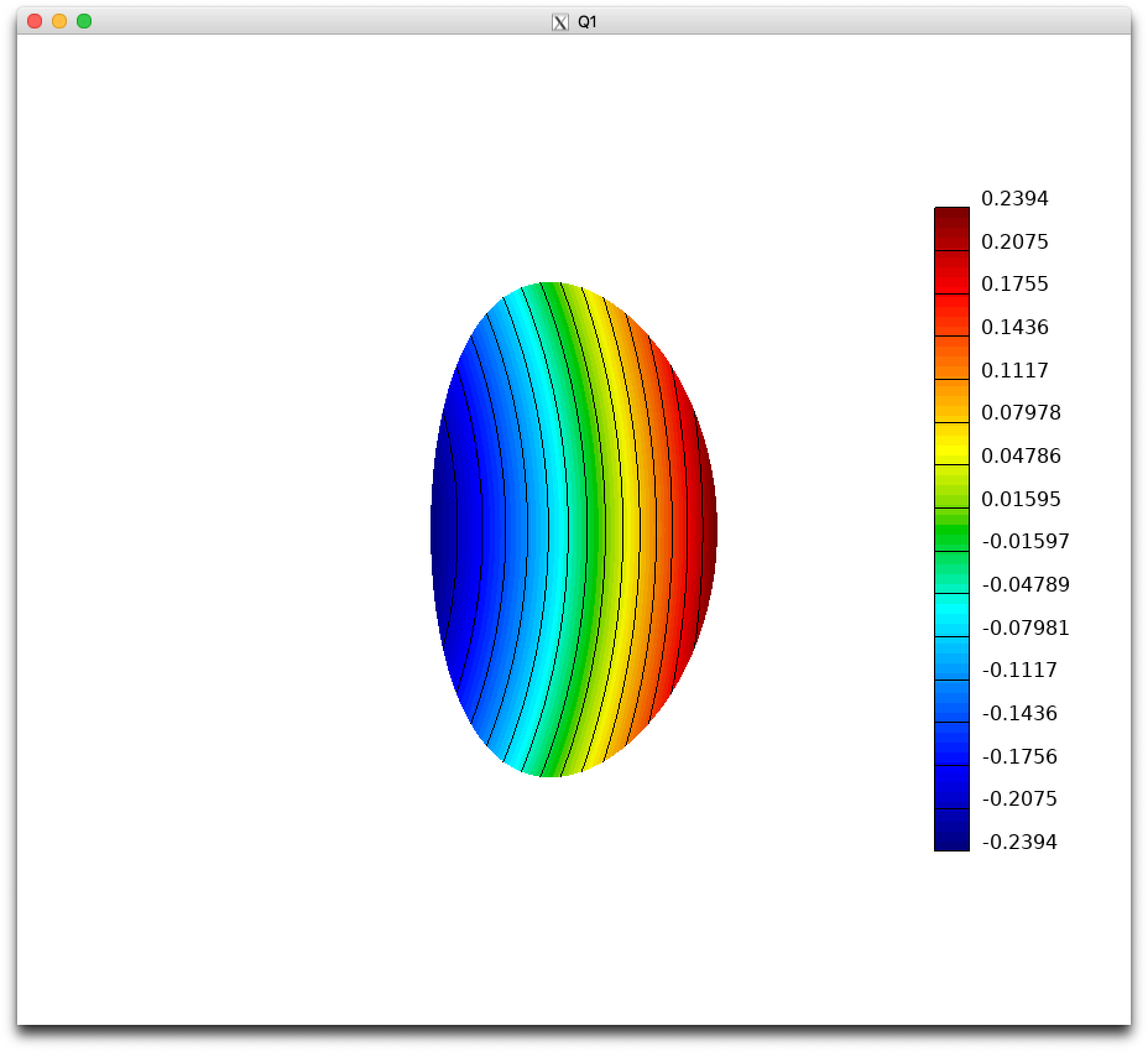}{ITER}{}{$\psi_r/r$}{-0.208}{0.239};
 \drawContour{xshift=10.5cm,yshift=6.75cm}{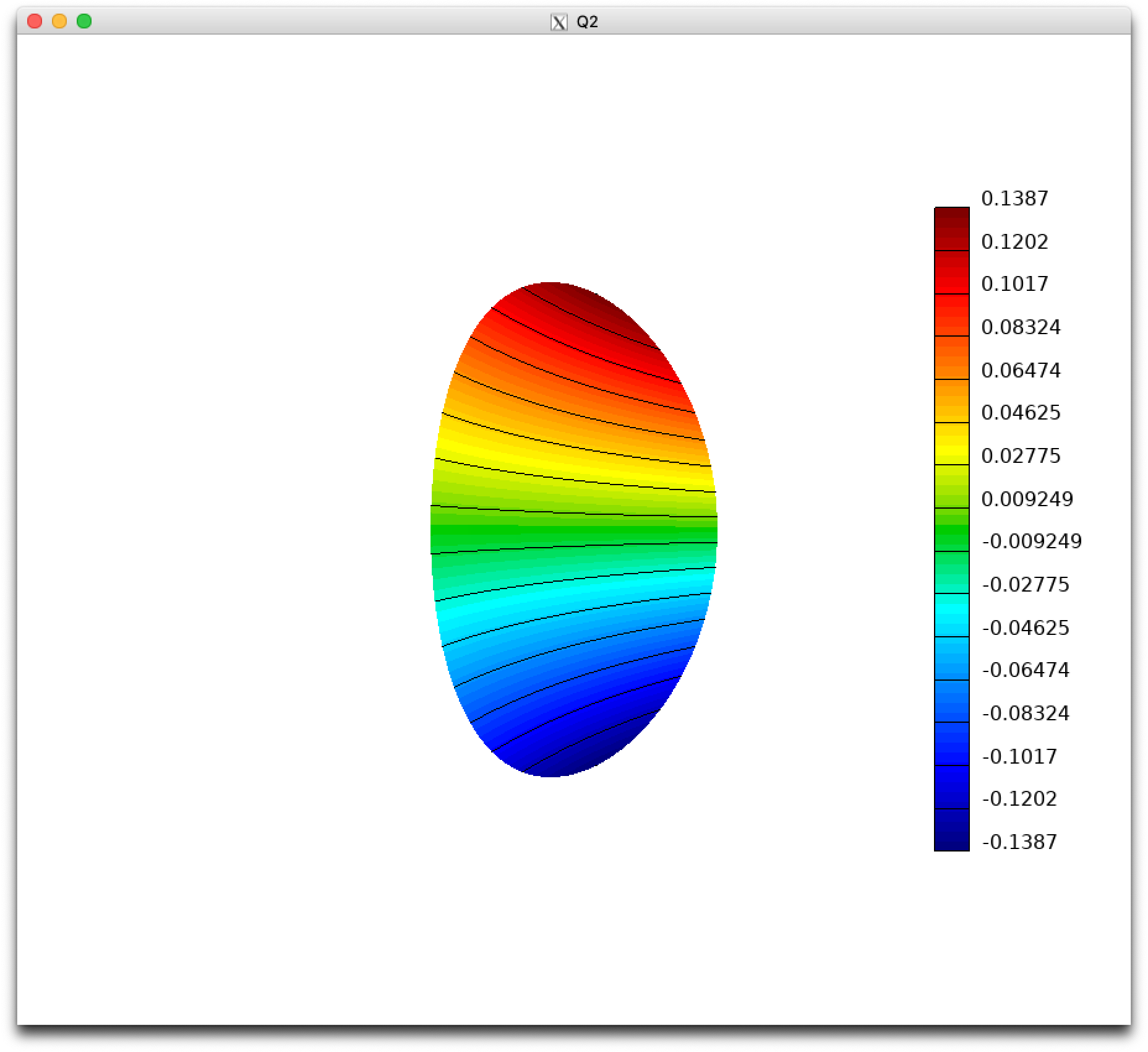}{ITER}{}{$\psi_z/r$}{-0.120}{0.139};
 \drawContour{xshift= 0.0cm,yshift=0.5cm}{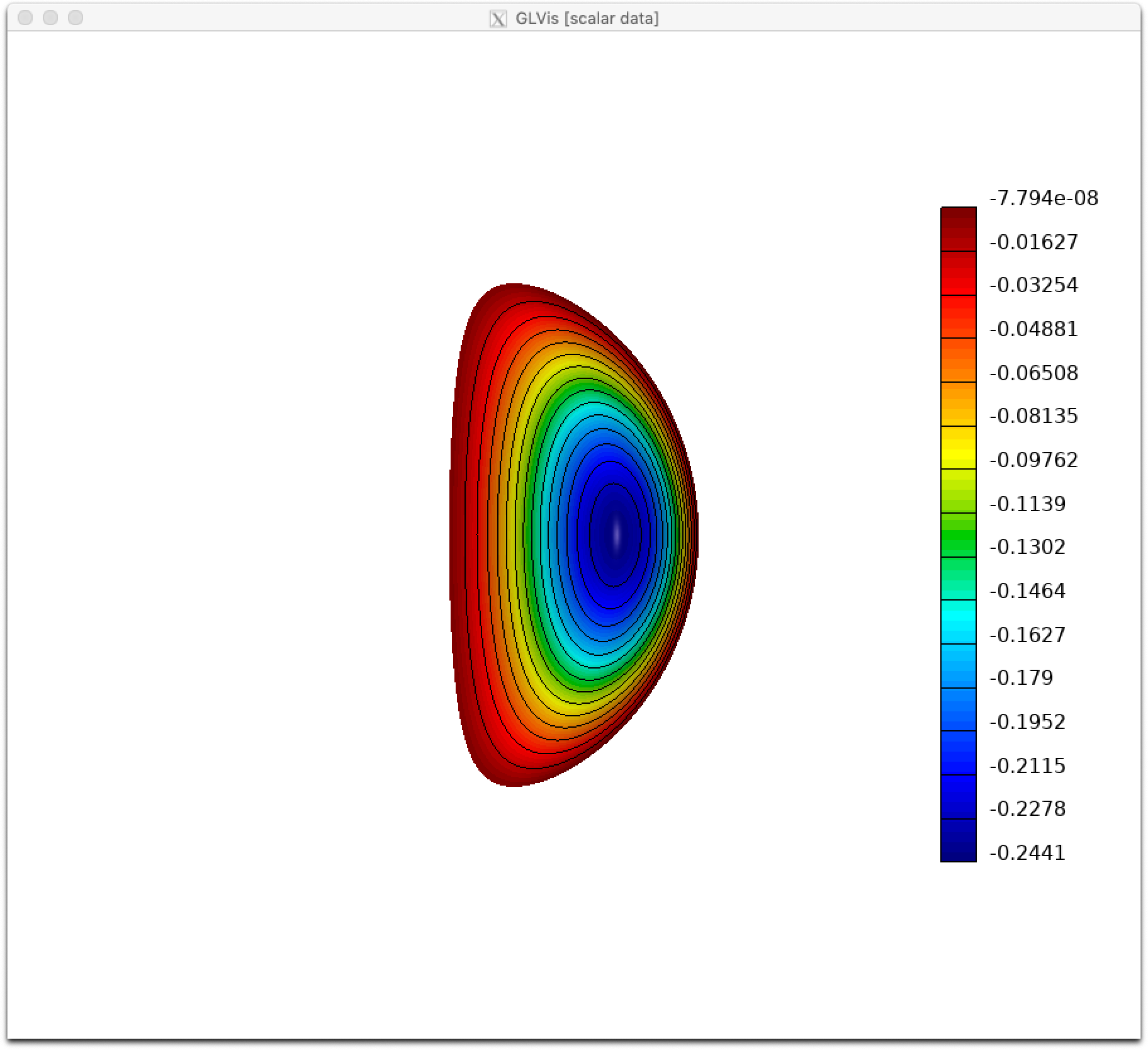}{NSTX}{}{$\psi$}{-0.227}{0};
 \drawContour{xshift=5.25cm,yshift=0.5cm}{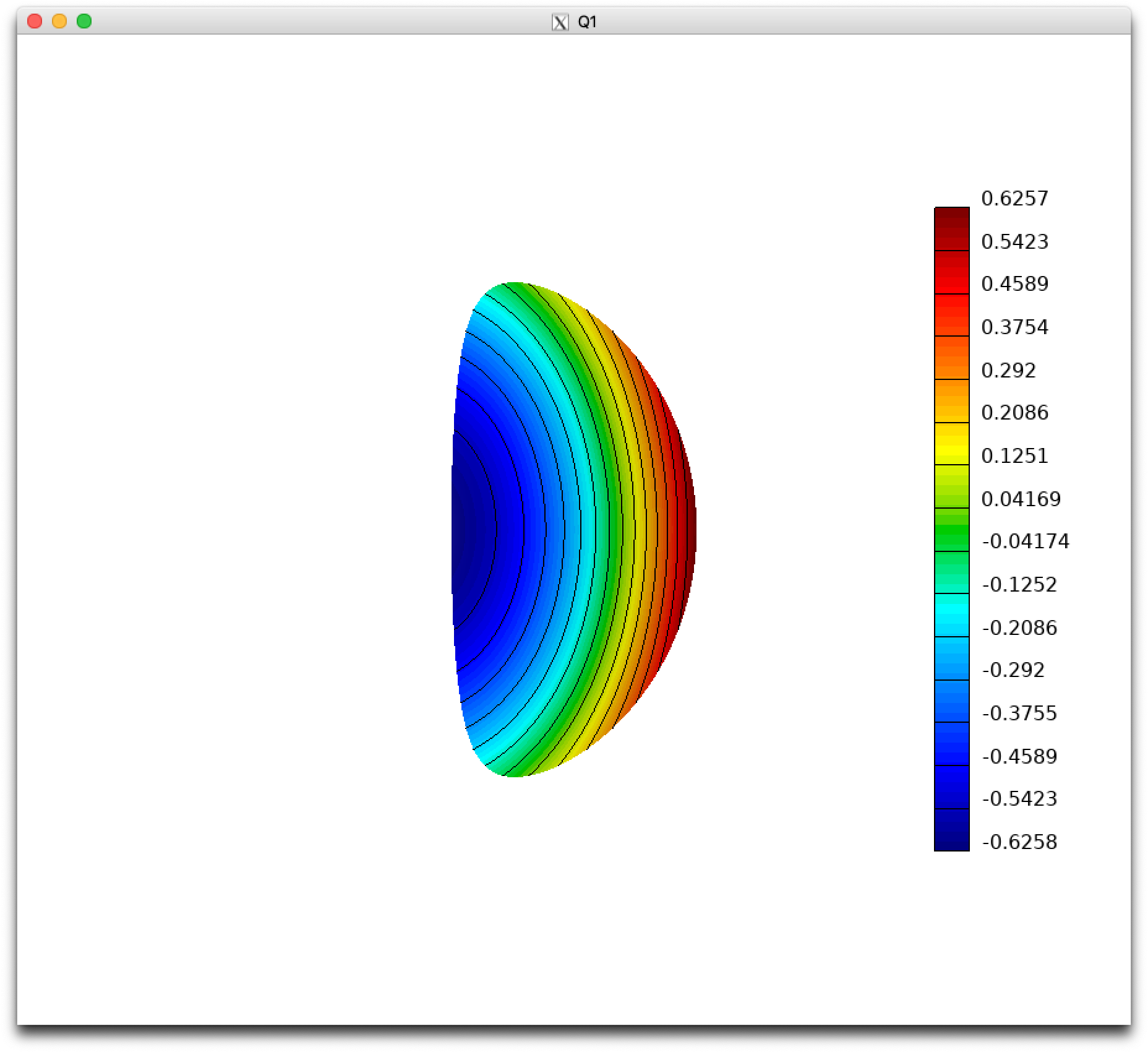}{NSTX}{}{$\psi_r/r$}{-0.542}{0.626};
 \drawContour{xshift=10.5cm,yshift=0.5cm}{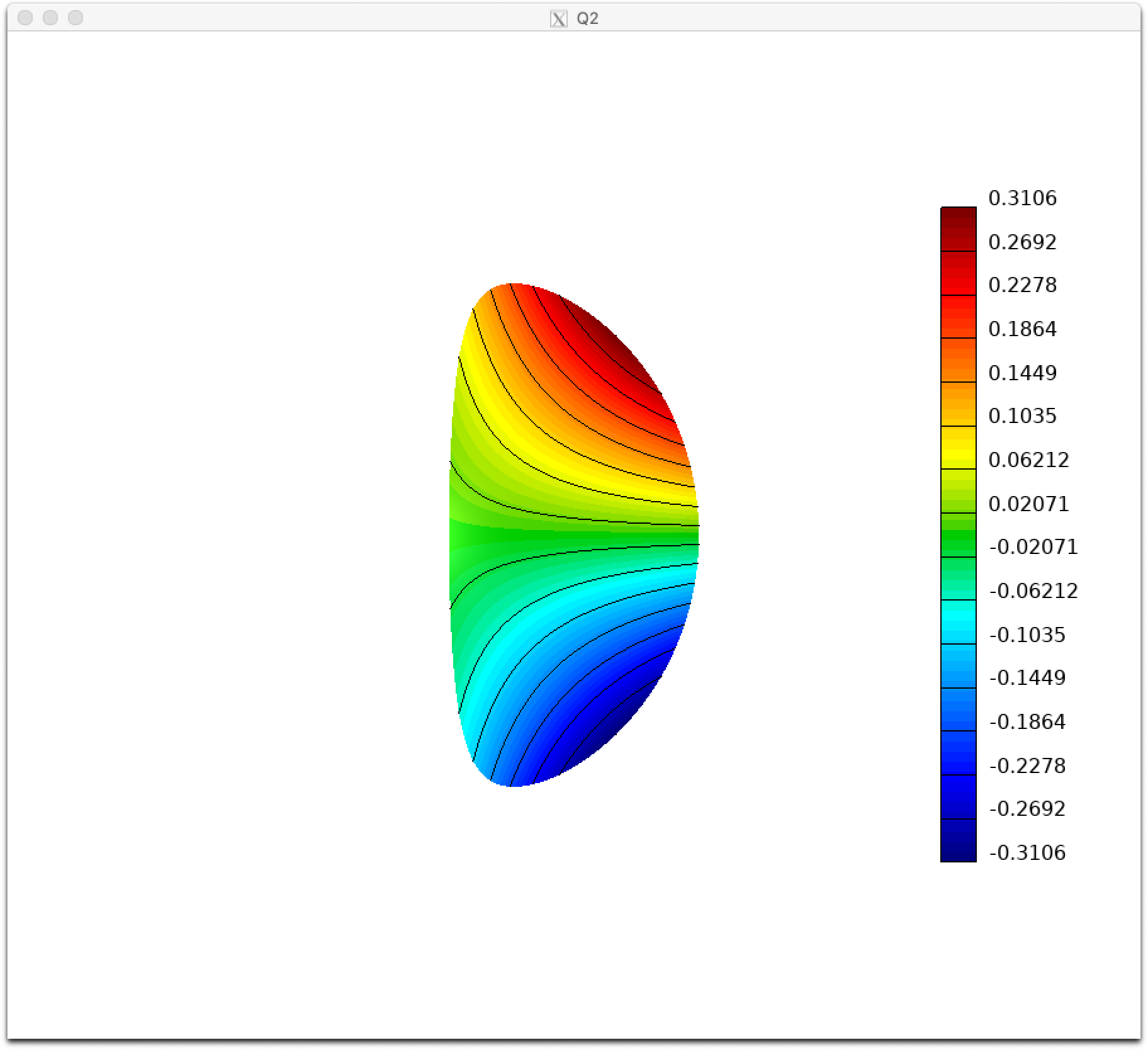}{NSTX}{}{$\psi_z/r$}{-0.269}{0.311};
\end{tikzpicture}
}
\end{center}
\caption{Numerical solutions and their derivatives for two linear accuracy tests in Section~\ref{sec:accuracy_linear}.  
The DPG scheme is used with the meshes of $\Gc^{(1)}$ and quadratic polynomials.
}
\label{fig:accuracy1}
\end{figure}
}

This example is further used to study the accuracy of the DPG schemes
as well as several commonly used finite element schemes.  We presented
the numerical solutions solved by two other schemes, the HDG and CG
methods.  To approximate $\qv=\wnabla\psi/r$, note that the ultraweak
DPG and HDG methods directly solve $\qv_h$ in the weak formulations,
while the CG method only solve $\psi_h$ and its $\qv_h$ is obtained by
solving $(\qv_h,\phi) = (\wnabla \psi_h/ r,\phi)$ {for arbitrary
  $\phi$ in the discrete test space.}  It is well known that such an
approach to compute $\qv_h$ will result in accuracy reduction by at
least one order.  On the other hand, the results of the HDG method do
not perform the post-processing reconstruction, which typically could
lift {one order of accuracy for $\psi$}.  To demonstrate the results,
as an example, quadratic polynomials are used for the trial space of
all three schemes.  The uniform refinement is performed during the
convergence study and the initial meshes are given in
Figure~\ref{fig:accuracy_mesh}.  $L^\infty$-errors of $\psi$ and $\qv$
are presented for different schemes in Table~\ref{tab:accuracy1-q}.
Here the $L^\infty$-error is defined as
\begin{align*}
L^{\infty}_{\textrm{error} }(u_h) = \max_{K\in \Omega_h} || u_h(x) - u_\textrm{exact}(x) ||_{\infty,K},
\end{align*}
where $||\cdot||_{\infty,K}$ is the standard $L^{\infty}$ norm on $K$.
The $L_\infty$-error for the vector $\qv$ is taken as the maximal error among all its components.
All three schemes achieve a third-order accuracy for $\psi$ as expected. The accuracy for $\qv$ of the DPG and HDG methods are third order, while the CG method only has second order accuracy. Compared with the HDG method, though with the same accuracy order, the error of the DPG method is smaller for both $\psi$ and $\qv$. 
We note that the HDG method is computationally more efficient than the DPG scheme for this case, as it results into a linear system of smaller size through a Schur complement.
The Schur complement for the matrix of $B^TG^{-1}B$ in the linear DPG scheme is not as obvious as the HDG scheme. 
It becomes even more challenging for the nonlinear problems we will consider later. 
Therefore, we do not investigate this direction in this work. {It is also important to \pzc{mention} that the post-processing reconstruction for the HDG scheme can not improve the accuracy of $\qv=\wnabla\psi/r$, though it lifts the accuracy order of $\psi$.}

The convergence results indicate the DPG scheme produce more accurate derivates than the other two schemes.

\begin{table}[htbp]
  \label{tab:accuracy1-q}%
  \centering
  \caption{$L^\infty$-errors and orders of  $\psi$ and $\qv=-\wnabla \psi/r$ for the linear problems in Section~\ref{sec:accuracy_linear}.   The DPG, HDG and CG schemes 
   are used with $P^2$ polynomials. It demonstrates the accuracy of the DPG scheme and also shows the performance of the
  different schemes when computing its derivative $\qv$. }
  \medskip
      \begin{tabular}{|c|c|c|c|c|c|c|c|c|c|c|c|c|c|}
    \hline
     \multicolumn{8}{|c|}{$L^\infty$-errors and orders of  $\psi$} \\ \hline 
\multirow{2}{*}{Test} &\multirow{2}{*} {Grid}&
       \multicolumn{2}{|c|}{ DPG}& \multicolumn{2}{|c|}{ HDG}&\multicolumn{2}{|c|}{ CG}\\
       \cline{3-8}
	  & 				& $\psi$  & order 			& $\psi$  & order		& $\psi$ & order \\ \hline
       	  &$\Gc^{(1)}$	&  2.877e-07   &-		&4.611e-07&-	   	&3.242e-06&- \\
 ITER &$\Gc^{(2)}$	&  3.609e-08  &2.99 		&6.035e-08&2.93  	&4.236e-07&2.94\\
 	  &$\Gc^{(4)}$	&  4.532e-09 & 2.99		&7.745e-09&2.96   	&5.335e-08&2.99\\
       	  &$\Gc^{(8)}$	&  5.682e-10 & 3.00 		&9.809e-10& 2.98 	&6.688e-09&3.00\\ \hline
       	  &$\Gc^{(1)}$	&  4.476e-06  &-		&6.238e-06&-	   	& 4.666e-05&- \\
NSTX &$\Gc^{(2)}$	&  6.062e-07  &2.88	 	&9.802e-07&2.67  	& 6.500e-06&2.84\\
 	  &$\Gc^{(4)}$	&  7.790e-08  &2.96		&1.288e-07&2.93   	& 8.502e-07&2.93\\
      	  &$\Gc^{(8)}$	&  9.874e-09  &2.98		&1.670e-08&2.95 	& 1.102e-07&2.95\\ \hline
    \end{tabular}%
    
  \bigskip
  
    \begin{tabular}{|c|c|c|c|c|c|c|c|c|c|c|c|c|c|}
        \hline
     \multicolumn{8}{|c|}{$L^\infty$-errors and orders of  $\qv$} \\ \hline 
\multirow{2}{*}{Test} &\multirow{2}{*} {Grid}&
       \multicolumn{2}{|c|}{ DPG}& \multicolumn{2}{|c|}{ HDG}&\multicolumn{2}{|c|}{ CG}\\
       \cline{3-8}
	  & 				& $\qv$  & order 			& $\qv$  & order		& $\qv$ & order \\ \hline
       	  &$\Gc^{(1)}$	&  6.742e-07   &-		&6.338e-06&-	   	&4.248e-04&- \\
 ITER &$\Gc^{(2)}$	&  8.696e-08  &2.95 		&8.162e-07&2.96  	&1.210e-04&1.81\\
 	  &$\Gc^{(4)}$	&  1.084e-08 & 3.00		&1.036e-08&2.98   	&3.129e-05&1.95\\
       	  &$\Gc^{(8)}$	&  1.369e-09 & 2.99 		&\pzc{1.322e-09}& 2.97 	&7.798e-06&2.00\\ \hline
       	  &$\Gc^{(1)}$	&  4.858e-05  &-		&1.088e-04&-	   	& 8.700e-03&- \\
NSTX &$\Gc^{(2)}$	&  6.384e-06  &2.93	 	&1.548e-05&2.81  	& 2.194e-03&1.99\\
 	  &$\Gc^{(4)}$	&  8.188e-07  &2.96		&2.039e-06&2.02   	& 5.468e-04&2.00\\
      	  &$\Gc^{(8)}$	&  1.037e-07  &2.98		&2.634e-07&2.05 	& 1.342e-04&2.03\\ \hline
    \end{tabular}%
\end{table}%

\subsubsection{Nonlinear test}
\label{sec:accuracy_nonlinear}
Next consider a manufactured solution given in \cite{sanchez2019hybridizable}. 
The solution satisfies 
\begin{align*}
\psi(r,z) =& \sin \left(k_r(r+r_0) \right) \, \cos(k_z z),
\end{align*}
with the source term given by
\begin{align*}
F(r,z,\psi) =& (k_r^2+k_z^2)\psi + \frac{k_r}{r} \cos( k_r(r+r_0) )\cos(k_z z) 
+r\Big[ \sin^2\left(k_r(r+r_0)\right) \cos^2(k_z z) \\ -& \psi^2 + \pzc{\exp\left(-\sin\big(k_r(r+r_0) \big)\cos(k_z z)\;\right)}-\exp(-\psi)\Big],
\end{align*}
where  the coefficients are $k_r=1.15\pi$, $k_z= 1.15$ and $r_0=-0.5$. 
We consider the same ITER geometry as the first linear test in the previous section. 
The computational domain and the initial mesh are identical to the ITER case given in Figure \ref{fig:accuracy_mesh}.
For this case, a non-homogenous Dirichlet boundary condition given by the exact solution is used throughout the computations.
The solutions are presented in Figure \ref{fig:accuracy2}, in which quadratic polynomials are used for the trial space.

 We further perform  a convergence study through a uniform refinement.
 This case is also used to demonstrate the capability of the arbitrary-order DPG scheme in our implementation.
To test the accuracy, $U_h^k$ with $k=1,2,3$ are considered for the discrete trial space. 
$L^\infty$ numerical errors and the corresponding orders are presented in Table~\ref{tab:accuracy2-q}. 
For the discrete trial space $U_h^k$ with $k=1,2$,
the optimal convergence of $(k+1)$-th order is observed for both the solutions and their derivatives. For $U_h^3$, though order reduction of $\psi$ occurs, optimal convergence for its derivatives is still observed.

 {%
\newcommand{\drawContour}[7]{%
\begin{scope}[#1]
\draw(0.0,-.3) node[anchor=south west,xshift=-4pt,yshift=+0pt] {\trimfig{picture/#2}{\figWidth}};
\draw(.5,5.) node[draw,fill=white,anchor=west,xshift=2pt,yshift=0pt] {\small #5};
\begin{scope}[xshift=2pt]
  \draw (\xcb,\ycb) node[anchor=south west,xshift=0.cm,yshift=.5cm,rotate=-90] {\trimfigcb{picture/test1_psi.png}{\cbWidth}{\cbHeight}};
  \draw (.8,0) node[anchor=north,xshift=+3pt,yshift=-2pt] {\scriptsize $#6$};
  \draw (4.4,0) node[anchor=north,xshift=+0pt,yshift=-2pt] {\scriptsize $#7$};
\end{scope}
\end{scope}
}
\newcommand{\cbWidth}{.3cm}
\newcommand{\cbHeight}{4cm}
\newcommand{\xcb}{.5cm}
\newcommand{\ycb}{-.2cm}
\setlength{\ycbTop}{\ycb+\cbHeight}
\setlength{\ycbMid}{\ycb+\cbHeight*\real{.5}}
\newcommand{\trimfigcb}[3]{\includegraphics[width=#2, height=#3, clip, trim=26.5cm 4cm 5cm 4cm]{#1}}
%
\newcommand{\figWidth}{5cm}
\newcommand{\yb}{6.75}
\newcommand{\trimfig}[2]{\trimw{#1}{#2}{.25}{.25}{.2}{.2}}
\begin{figure}[htb]
\begin{center}
\resizebox{13cm}{!}{
\begin{tikzpicture}[scale=1]
  \useasboundingbox (0.0,.0) rectangle (16.,6);  
 \drawContour{xshift= 0.0cm,yshift=0cm}{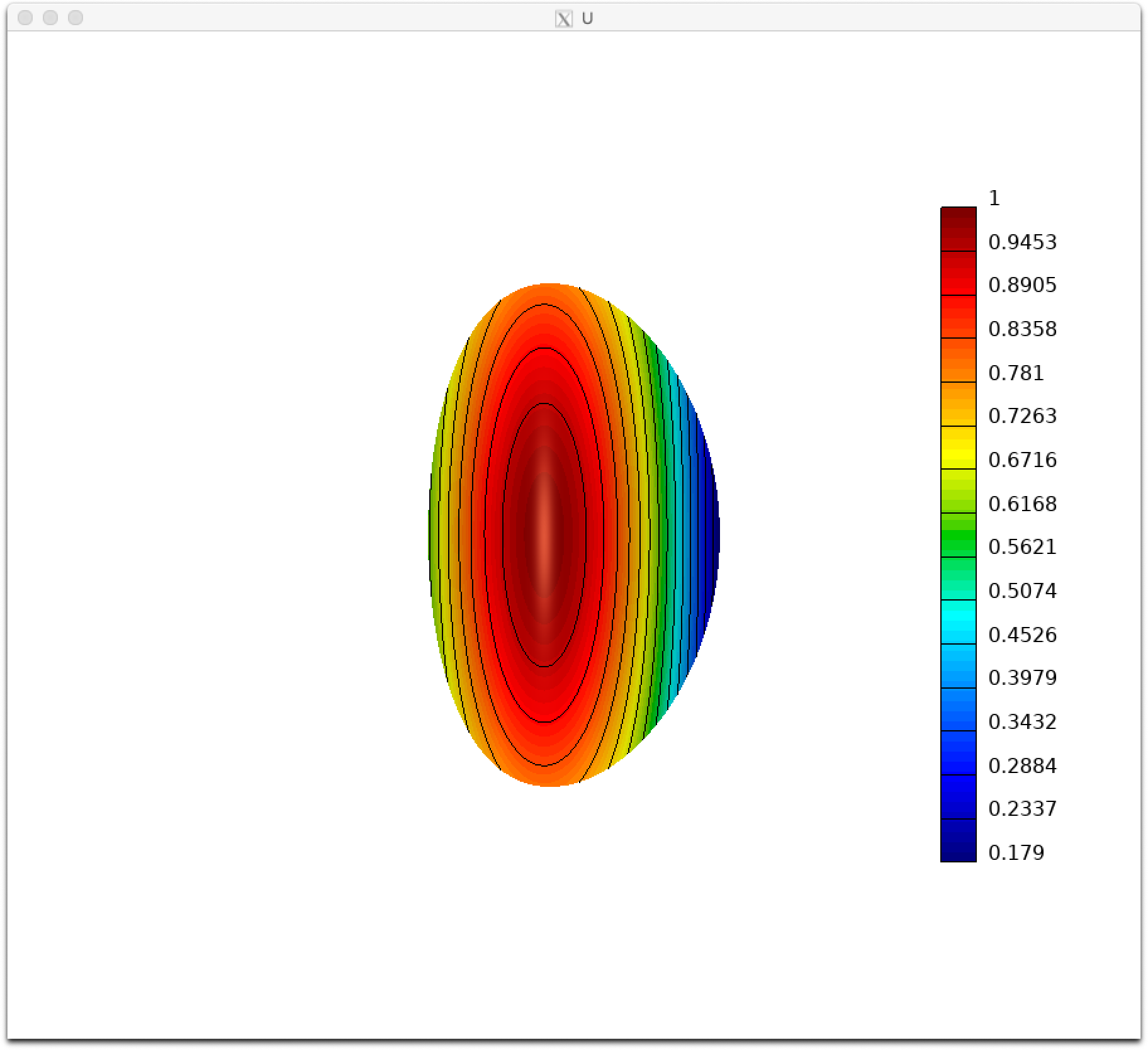}{}{}{$\psi$}{0.234}{1};
 \drawContour{xshift=5.25cm,yshift=0cm}{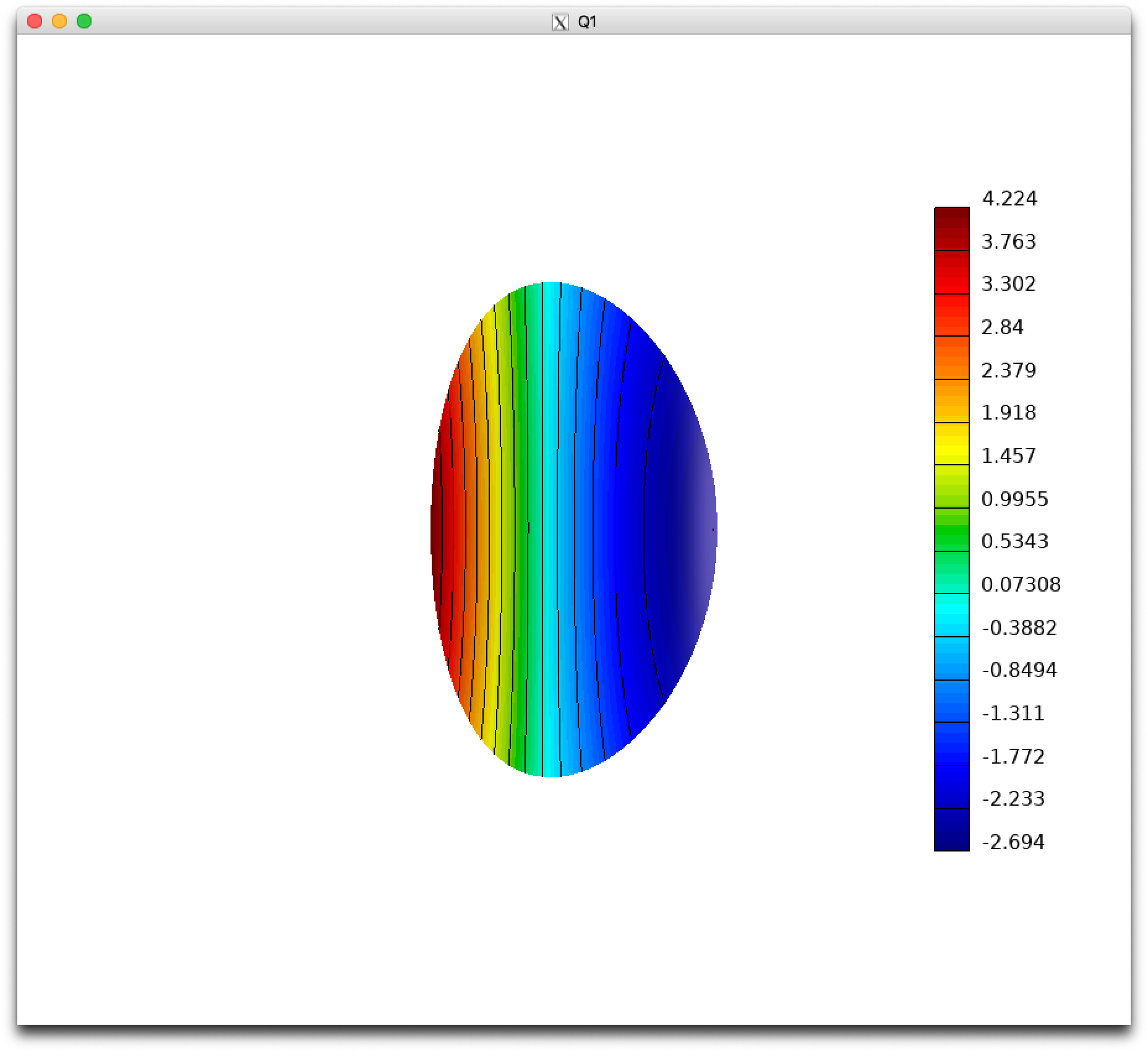}{}{}{$\psi_r/r$}{-2.23}{4.22};
 \drawContour{xshift=10.5cm,yshift=0cm}{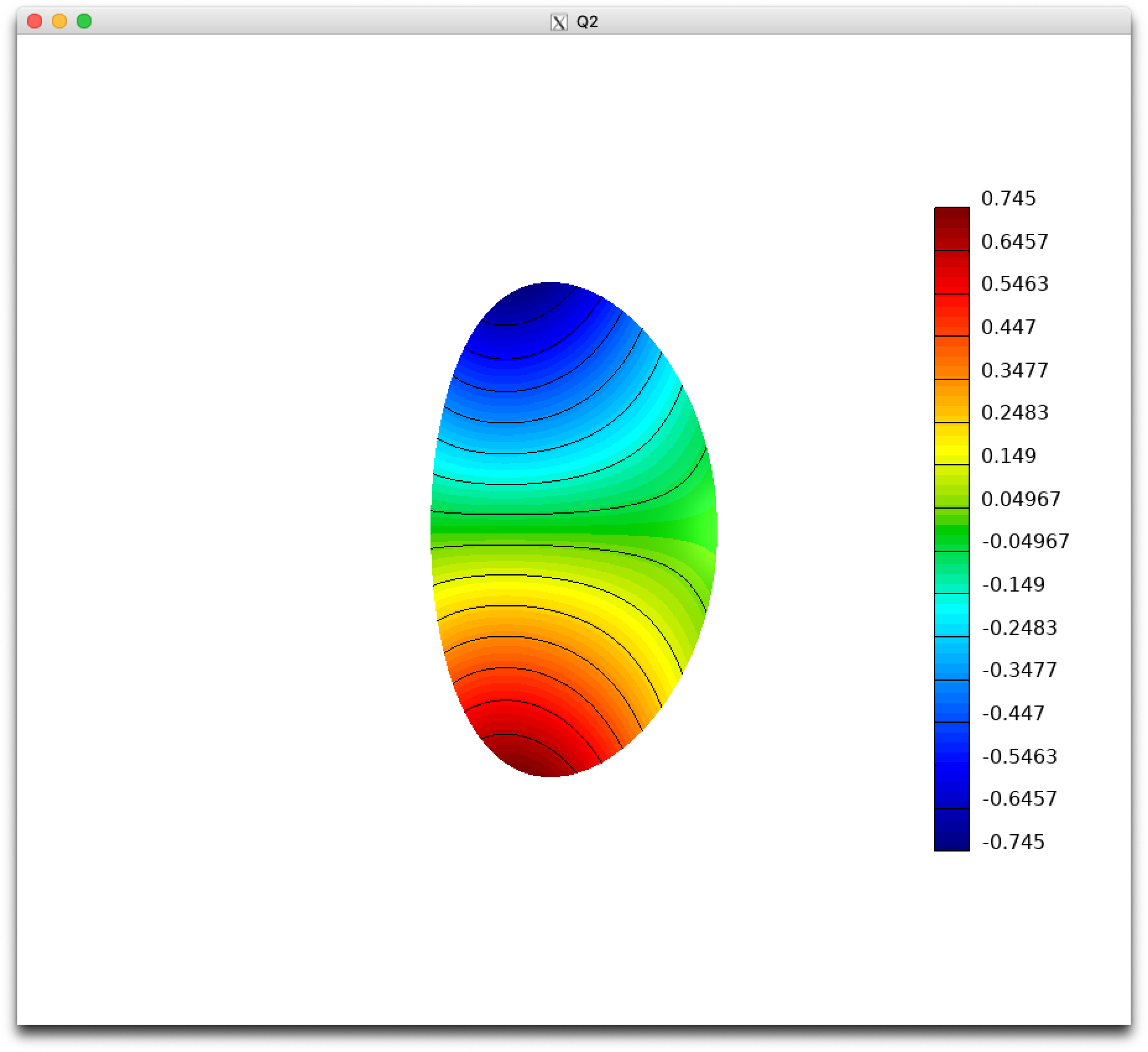}{}{}{$\psi_z/r$}{-0.646}{0.745};
\end{tikzpicture}
}
\end{center}
\caption{Numerical solutions and their derivatives for the nonlinear accuracy problem in Section~\ref{sec:accuracy_nonlinear}.  
The DPG scheme is used with the meshes of $\Gc^{(1)}$ and quadratic polynomials.
}
\label{fig:accuracy2}
\end{figure}
}

\begin{table}[htbp]
  \label{tab:accuracy2-q}%
  \centering
  \caption{$L^\infty$ errors and orders of   $\psi$ and $\qv=-\wnabla \psi/r$ for the nonlinear problem in Section~\ref{sec:accuracy_nonlinear}.
  $P^1$, $P^2$ and $P^3$ polynomial spaces are used for the DPG scheme.
  }
  \medskip
      \begin{tabular}{|c|c|c|c|c|c|c|c|c|c|c|c|c|c|}
        \hline
     \multicolumn{7}{|c|}{$L^\infty$-errors and orders of  $\psi$} \\ \hline 
\multirow{1}{*} {Grid}&
       \multicolumn{2}{|c|}{ $P^1$}& \multicolumn{2}{|c|}{ $P^2$}&\multicolumn{2}{|c|}{ $P^3$}\\ \hline
	$\Gc^{(1)}$	&  1.074e-03   &-		& 2.892e-05&-	   	& 2.294e-07&- \\
	$\Gc^{(2)}$	&  2.648e-04  & 2.02 	& 2.604e-06& 3.47  	& 1.314e-08&4.13\\
	$\Gc^{(4)}$	&  6.572e-05  & 2.01		& 2.448e-07& 3.41   	& 5.446e-09&1.27\\
	$\Gc^{(8)}$	&  1.639e-05 &  2.00 	& 2.451e-08& 3.32 	& 1.720e-08&-1.66\\ \hline
    \end{tabular}%

\bigskip
  
    \begin{tabular}{|c|c|c|c|c|c|c|c|c|c|c|c|c|c|}
        \hline
     \multicolumn{7}{|c|}{$L^\infty$-errors and orders of  $\qv$} \\ \hline 
\multirow{1}{*} {Grid}&
       \multicolumn{2}{|c|}{ $P^1$}& \multicolumn{2}{|c|}{ $P^2$}&\multicolumn{2}{|c|}{ $P^3$}\\ \hline
	$\Gc^{(1)}$	&  2.226e-03   &-		& 4.338e-05&-	   	& 7.934e-07&- \\
	$\Gc^{(2)}$	& 5.685e-04  & 1.97 		& 4.697e-06& 3.21  	& 5.174e-08&3.94\\
	$\Gc^{(4)}$	& 1.433e-04  & 1.99		& 5.729e-07& 3.04   	& 3.321e-09&3.96\\
	$\Gc^{(8)}$	&  3.598e-05 & 1.99 		& 7.978e-08& 2.84 	& 2.385e-10&3.80\\ \hline
    \end{tabular}%
\end{table}%

\subsection{Nonlinear tests}
\label{sec:nonlinear_test}
Next we consider a D-shaped geometry taken from~\cite{sanchez2019hybridizable}. 
The boundary of computational domain $\omega$ is determined by
$$
r(s) = 1+\eps \cos(s+\arcsin(\delta \sin(s) ), \quad z(s) = \eps\kappa \sin(s),
$$
where $s\in[0,2\pi]$, $\eps=0.32$, $\delta=0.33$ and $\kappa=1.7$. The source term is 
\begin{align}
F(r,z,\psi) = r^2\Big[1-\frac{1}{2}(1-\psi^2)^2 \Big].
\end{align}
A homogenous Dirichlet boundary condition is used in the simulation.
Numerical results with trial space $U_h^2$ are shown in Figure \ref{fig:DShape}.
This case demonstrates the capability to handle different geometry and {nonlinear source} in our DPG implementation.

  {%
\newcommand{\drawContour}[7]{%
\begin{scope}[#1]
\draw(0.0,-.3) node[anchor=south west,xshift=-4pt,yshift=+0pt] {\trimfig{picture/#2}{\figWidth}};
\draw(.5,5.) node[draw,fill=white,anchor=west,xshift=2pt,yshift=0pt] {\small #5};
\begin{scope}[xshift=2pt]
  \draw (\xcb,\ycb) node[anchor=south west,xshift=0.cm,yshift=.5cm,rotate=-90] {\trimfigcb{picture/test1_psi.png}{\cbWidth}{\cbHeight}};
  \draw (.8,0) node[anchor=north,xshift=+3pt,yshift=-2pt] {\scriptsize $#6$};
  \draw (4.4,0) node[anchor=north,xshift=+0pt,yshift=-2pt] {\scriptsize $#7$};
\end{scope}
\end{scope}
}
\newcommand{\cbWidth}{.3cm}
\newcommand{\cbHeight}{4cm}
\newcommand{\xcb}{.5cm}
\newcommand{\ycb}{-.2cm}
\setlength{\ycbTop}{\ycb+\cbHeight}
\setlength{\ycbMid}{\ycb+\cbHeight*\real{.5}}
\newcommand{\trimfigcb}[3]{\includegraphics[width=#2, height=#3, clip, trim=26.5cm 4cm 5cm 4cm]{#1}}
%
\newcommand{\figWidth}{5cm}
\newcommand{\yb}{6.75}
\newcommand{\trimfig}[2]{\trimw{#1}{#2}{.25}{.25}{.2}{.2}}
\begin{figure}[htb]
\begin{center}
\resizebox{13cm}{!}{
\begin{tikzpicture}[scale=1]
  \useasboundingbox (0.0,.0) rectangle (16.,6);  
 \drawContour{xshift= 0.0cm,yshift=0cm}{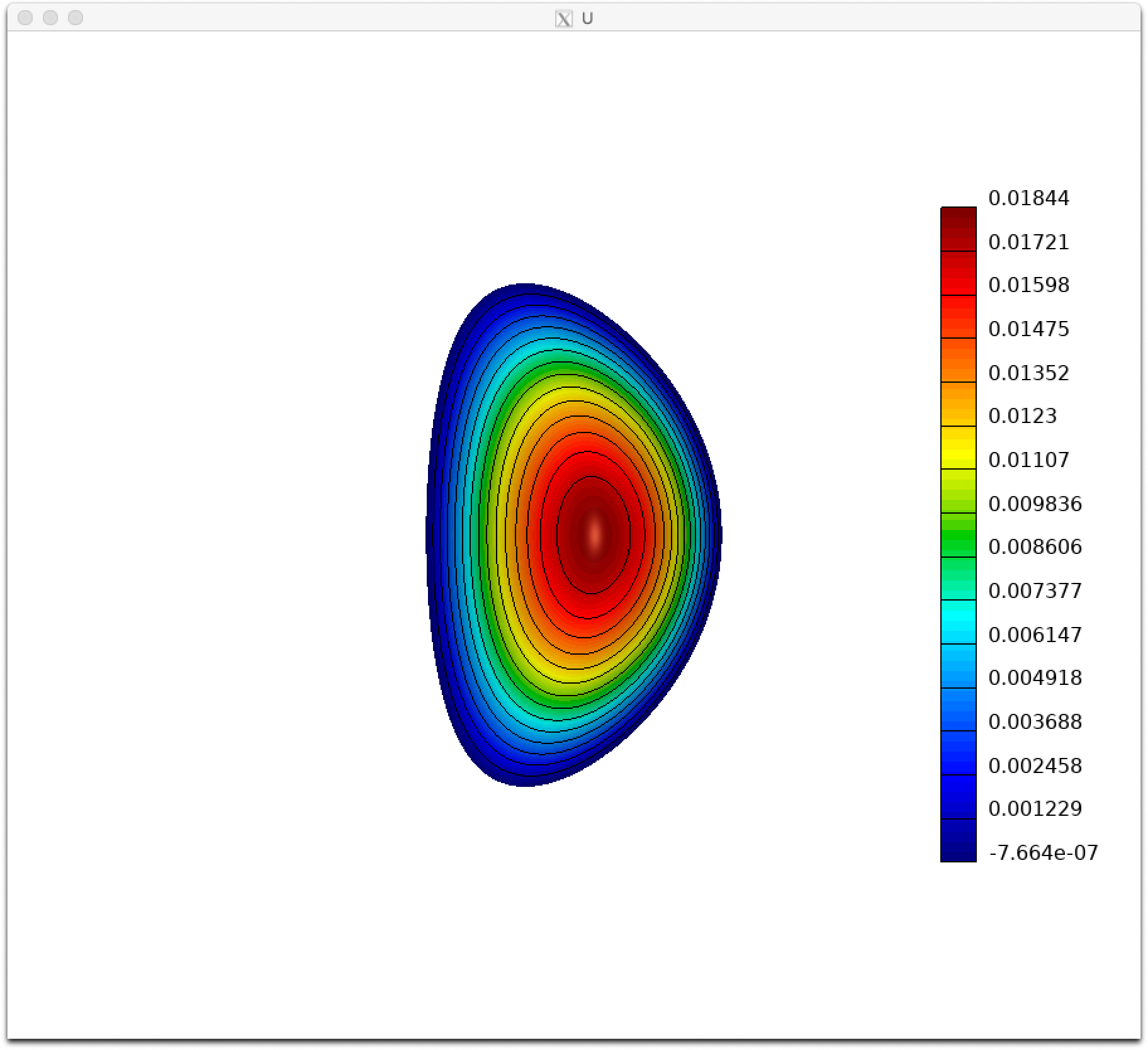}{}{}{$\psi$}{0.00123}{0.0184};
 \drawContour{xshift=5.25cm,yshift=0cm}{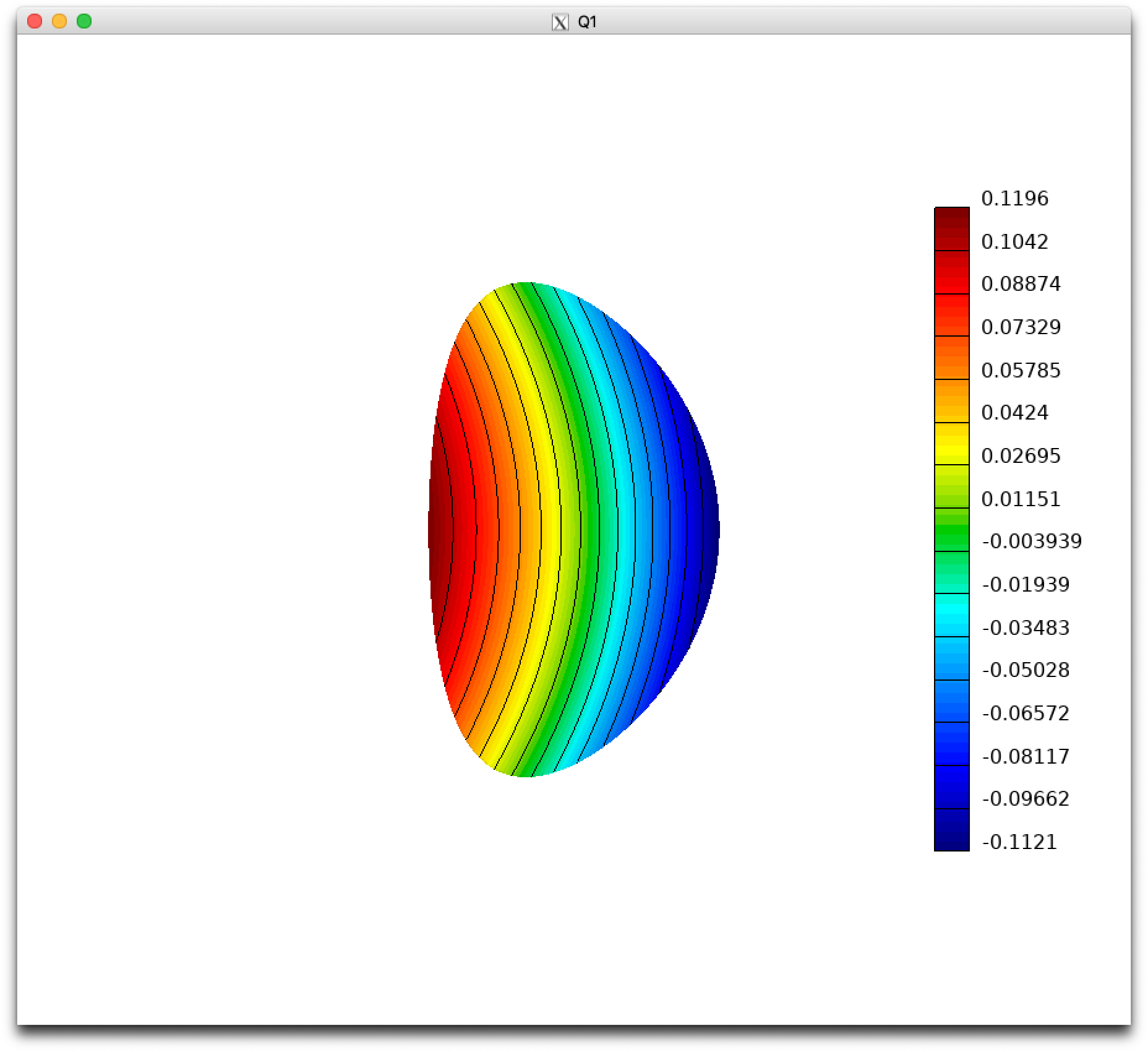}{}{}{$\psi_r/r$}{-0.0966}{0.120};
 \drawContour{xshift=10.5cm,yshift=0cm}{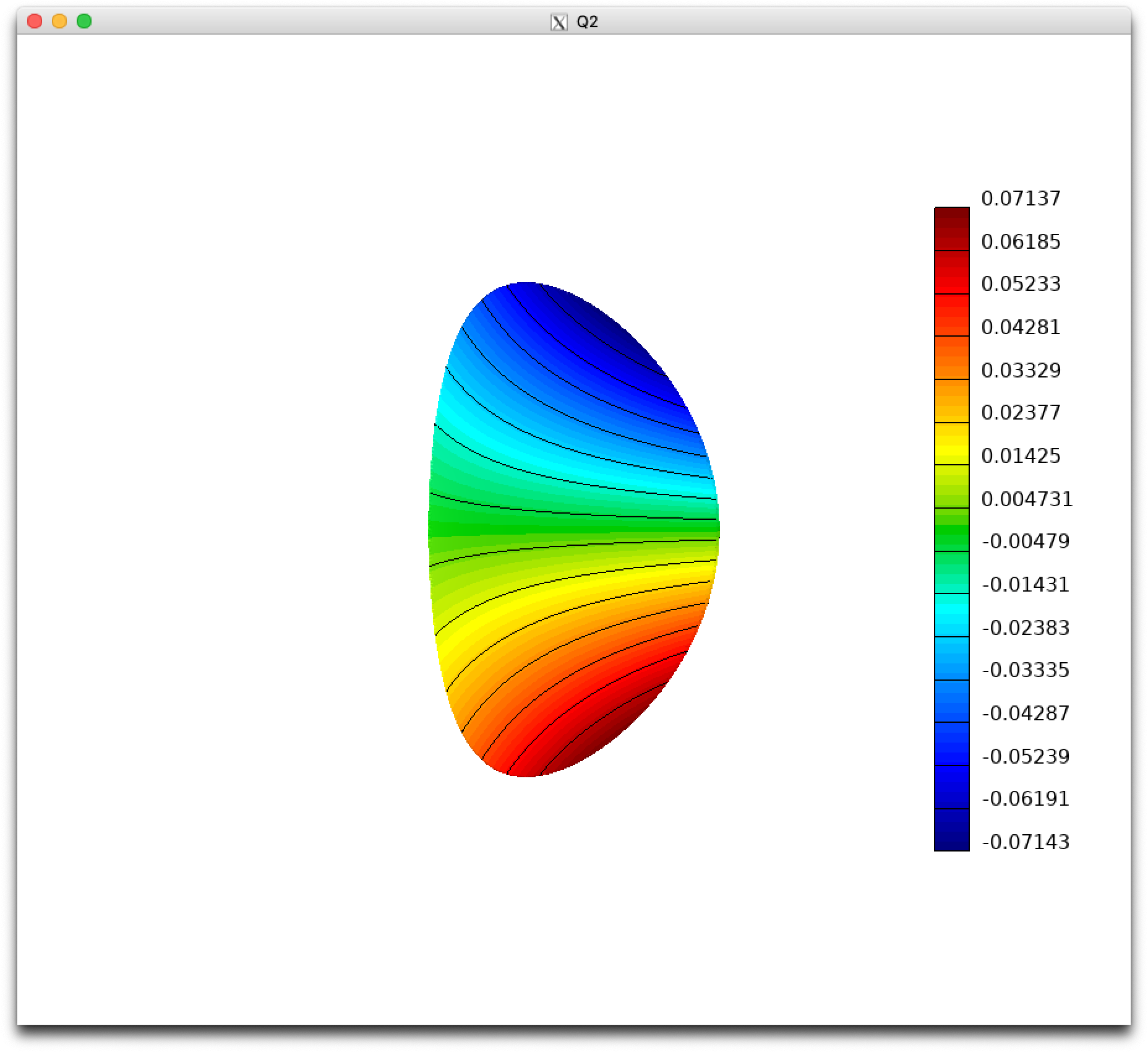}{}{}{$\psi_z/r$}{-0.0619}{0.0714};
\end{tikzpicture}
}
\end{center}
\caption{Numerical solutions for the case of the D-shaped geometry in Section~\ref{sec:nonlinear_test}.}
\label{fig:DShape}
\end{figure}
}

\subsection{Adaptivity}
\label{sec:amr_test}
In the final part of the numerical example, we focus on the performance of the DPG scheme when the AMR strategy is applied.
Two tests are presented with one involving a rectangular geometry and one involving D-shaped geometry.
We set $\theta_{\textrm{max} }=0.025$ and $\theta_{\textrm{global} }=0.025$. 
The choice of $\textrm{atol}_{\textrm{amr} }$ and the discrete trial space $U_h^k$ will be specified in each test.

{
\newcommand{\figWidth}{5.0cm}
\newcommand{\figWidthb}{4.9cm}
\newcommand{\trimfig}[2]{\trimh{#1}{#2}{.25}{.2}{.2}{.2}}
\newcommand{\figWidtha}{3.6cm}
\newcommand{\trimfiga}[2]{\trimh{#1}{#2}{.0}{.0}{.0}{.0}}
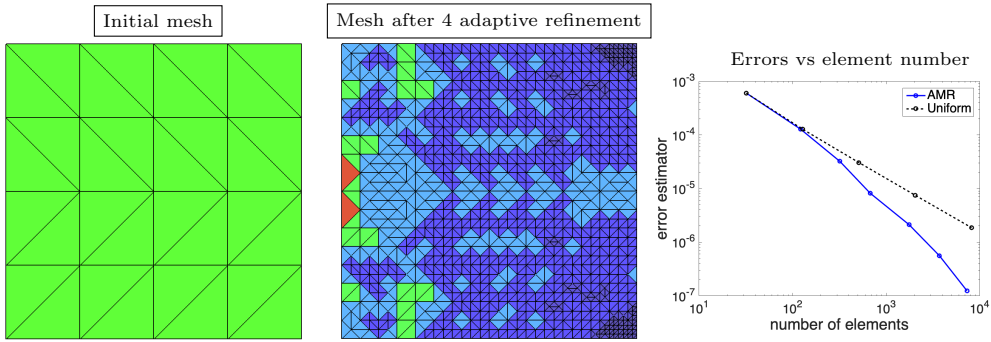
\begin{figure}[htb]
\begin{center}
\begin{tikzpicture}[scale=1]
  \useasboundingbox (0,.5) rectangle (12.,5.);  
  \draw(-1,-.3) node[anchor=south west,xshift=-4pt,yshift=+0pt] {\trimfig{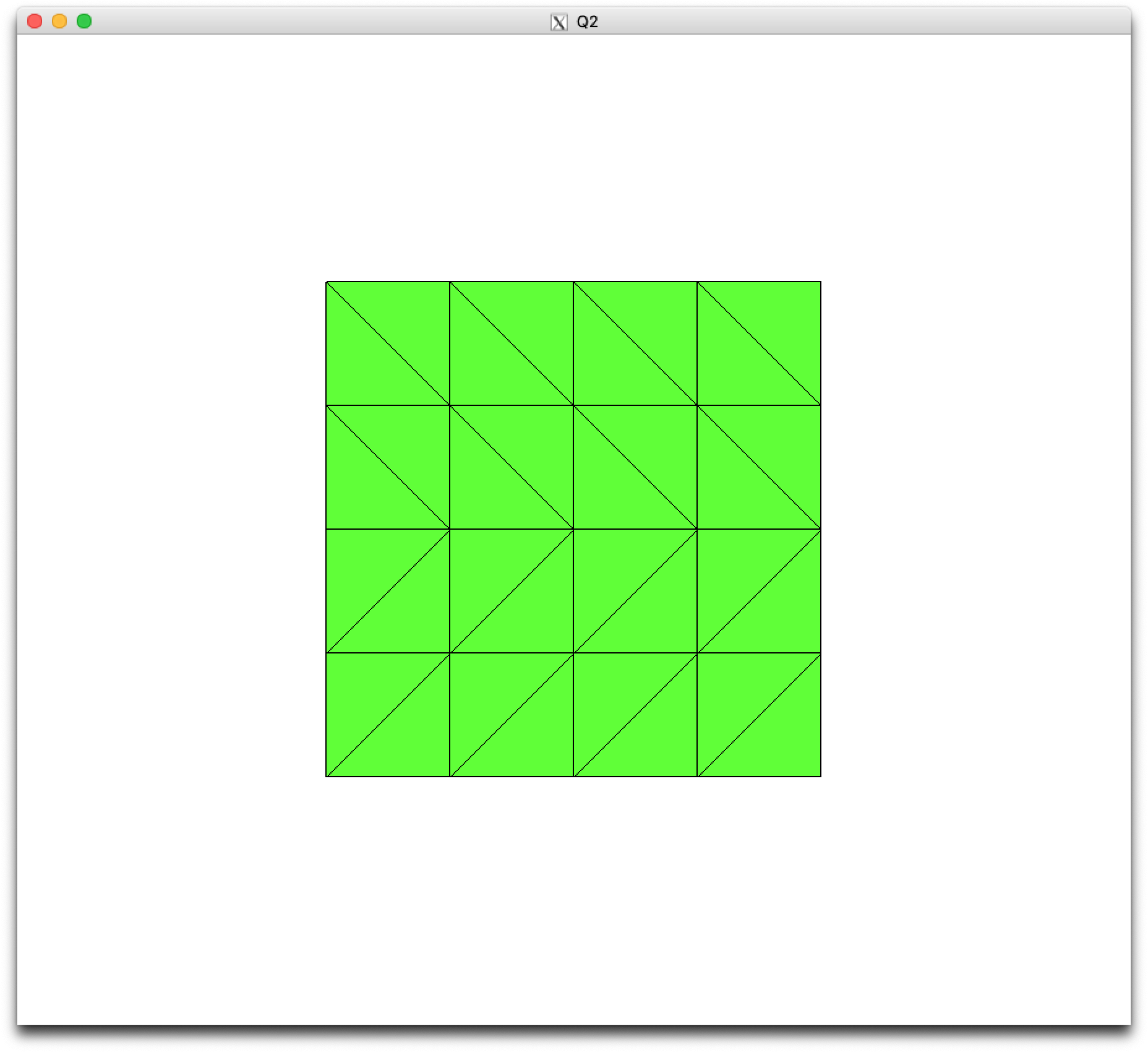}{\figWidth}};
  \draw(3.5,-.2) node[anchor=south west,xshift=-4pt,yshift=+0pt] {\trimfig{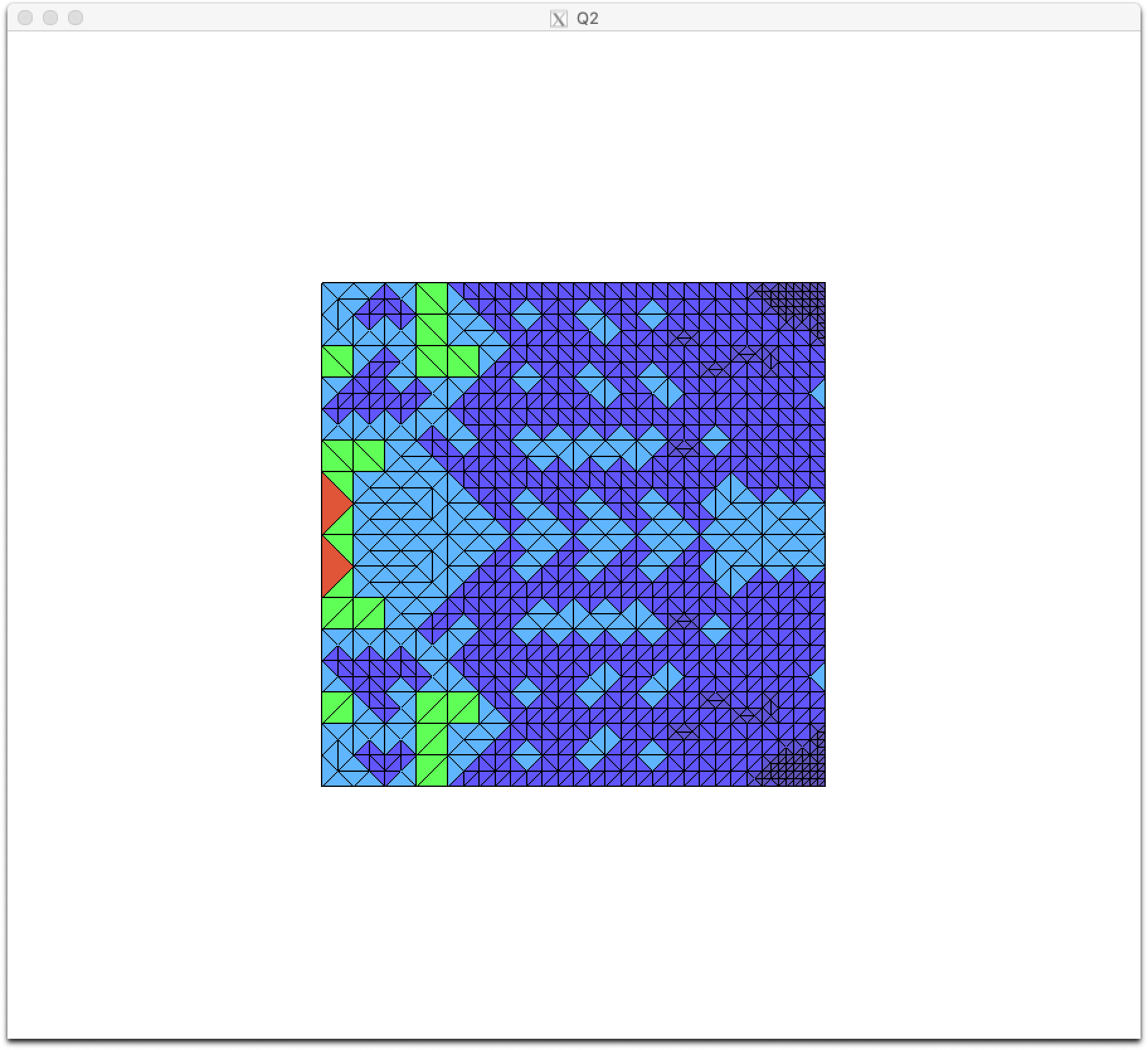}{\figWidthb}};
  \draw(7.9, .3) node[anchor=south west,xshift=-4pt,yshift=+0pt] {\trimfiga{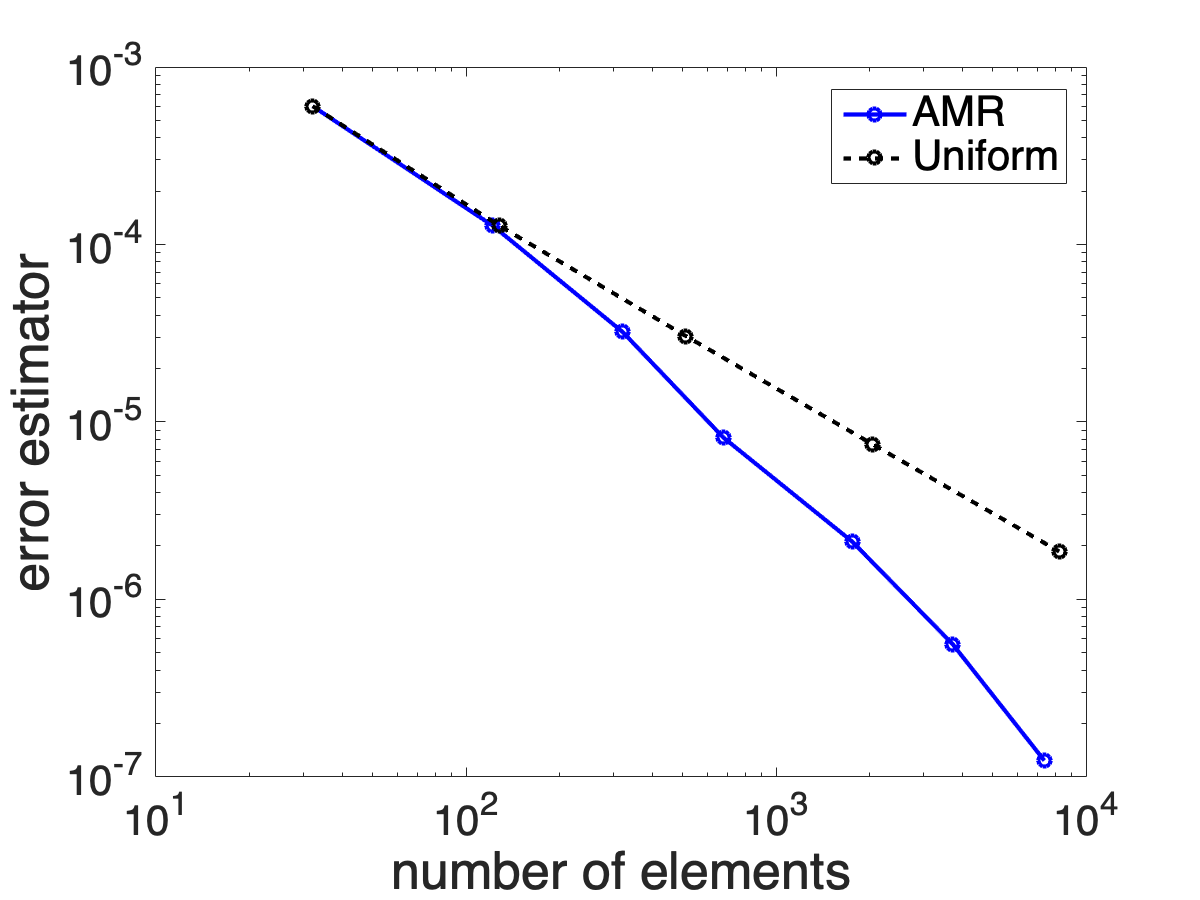}{\figWidtha}};
  \draw(1.7,4.6) node[draw,fill=white] {\scriptsize Initial mesh};
   \draw(6.1,4.6) node[draw,fill=white] {\scriptsize Mesh after 4 adaptive refinement};
    \draw(10.9,4.1) node {\scriptsize Errors vs element number};
    
%
\end{tikzpicture}
\end{center}
\caption{Left: initial mesh for the test in Section \ref{sec:amr_test1}. Middle: mesh after four adaptive refinement iterations. 
Right: errors verse total number of elements for both uniform and adaptive refinements. }
\label{fig:amr1_mesh}
\end{figure}
}

\subsubsection{Rectangular geometry}
\label{sec:amr_test1}
First consider a rectangle computational domain of $\Omega=[0.1,1.6]\times[-0.75,0.75]$ with a homogenous boundary condition of $\psi=0.25$.
The source term is  nonlinear given by
$$F(r,\psi) = 2r^2\psi\Big[ c_2(1-\exp(-\psi^2/\sigma^2)+\frac{1}{\sigma^2}(c_1+c_2\psi^2)\exp(-\psi^2/\sigma^2) \Big],$$ 
where $\sigma^2=0.005$,
$c_1=0.8$ and $c_2=0.2$. The reason for first considering the rectangular geometry is to eliminate the impact of the boundary condition and the extra complexity associated with a complex geometry
 when performing AMR.
The focus of this case is therefore on the improvement of numerical errors of AMR over a uniform mesh refinement.

For this test, we set $\textrm{atol}_{\textrm{amr} }=$1e-8 and use quadratic polynomial space $U_h^2$ as the discrete trial space. The initial mesh is presented in Figure~\ref{fig:amr1_mesh}.
Note that the initial mesh is chosen carefully so that it retains a symmetry along $y=0$, which is consistent with the exact solution.
In our solver, we perform one refinement check after each nonlinear solve, and the updated adaptive mesh will be used again in the next nonlinear solve.
An iteration procedure (which is referred to as AMR iterations in this work) is performed until a targeting error goal is achieved or the total number of elements is larger than an upper bound.

The mesh after four AMR iterations is presented in Figure \ref{fig:amr1_mesh}. 
It is accompanied by the numerical solutions on the initial and adaptive meshes presented in Figure \ref{fig:amr1}.
The color in the adaptive mesh of Figure \ref{fig:amr1_mesh} indicates different refinement levels of AMR.
We first find that the adaptive mesh remains symmetric throughout the AMR iterations. 
Note that $\qv=\wnabla\psi/r$ has sharp features near the top-right corner and the bottom-right corner, which are not well-resolved on the initial mesh. 
After four AMR iterations, these features are better captured and most refinements are performed near them.  
Overall, the computational efforts are focused on the region where interesting physics happens during the AMR iteration.
All those facts indicate our AMR strategy is effective and efficient.

As a further verification of the AMR strategy, we also compare the
numerical errors under uniform and adaptive refinements. {Here the
  numerical error is measured by the value of our error estimator
  defined in~\eqref{eqn:estimator_total}.}  The convergence histories
for both AMR and uniform refinement are also presented in Figure
\ref{fig:amr1_mesh}.  The numerical errors are presented under
increasing total number of elements.  Note that the errors on AMR are
much smaller than errors on uniform meshes when the total number of
elements (degree of freedom) are the same.  Another important
observations is that as the meshes are refined, the numerical errors
on AMR decay quicker than the errors on uniform meshes.  Therefore,
the convergence histories quantitatively confirm the efficiency of the
AMR approach in this work.


 {%
\newcommand{\drawContour}[7]{%
\begin{scope}[#1]
\draw(0.0,-.3) node[anchor=south west,xshift=-4pt,yshift=+0pt] {\trimfig{picture/amr/#2}{\figWidth}};
\draw(1.5,5.) node[draw,fill=white,anchor=west,xshift=2pt,yshift=0pt] {\small #3};
\draw(.5,5.) node[draw,fill=white,anchor=west,xshift=2pt,yshift=0pt] {\small #5};
\begin{scope}[xshift=2pt]
  \draw (\xcb,\ycb) node[anchor=south west,xshift=0.cm,yshift=.5cm,rotate=-90] {\trimfigcb{picture/test1_psi.png}{\cbWidth}{\cbHeight}};
  \draw (.8,0) node[anchor=north,xshift=+3pt,yshift=-2pt] {\scriptsize $#6$};
  \draw (4.6,0) node[anchor=north,xshift=+0pt,yshift=-2pt] {\scriptsize $#7$};
\end{scope}
\end{scope}
}
\newcommand{\cbWidth}{.3cm}
\newcommand{\cbHeight}{4.5cm}
\newcommand{\xcb}{.5cm}
\newcommand{\ycb}{-.2cm}
\setlength{\ycbTop}{\ycb+\cbHeight}
\setlength{\ycbMid}{\ycb+\cbHeight*\real{.5}}
\newcommand{\trimfigcb}[3]{\includegraphics[width=#2, height=#3, clip, trim=26.5cm 4cm 5cm 4cm]{#1}}
%
\newcommand{\figWidth}{5cm}
\newcommand{\yb}{6.75}
\newcommand{\trimfig}[2]{\trimw{#1}{#2}{.25}{.25}{.2}{.2}}
\begin{figure}[htb]
\begin{center}
\resizebox{13cm}{!}{
\begin{tikzpicture}[scale=1]
  \useasboundingbox (0.0,.5) rectangle (16.,12);  
 \drawContour{xshift= 0.0cm,yshift=6.75cm}{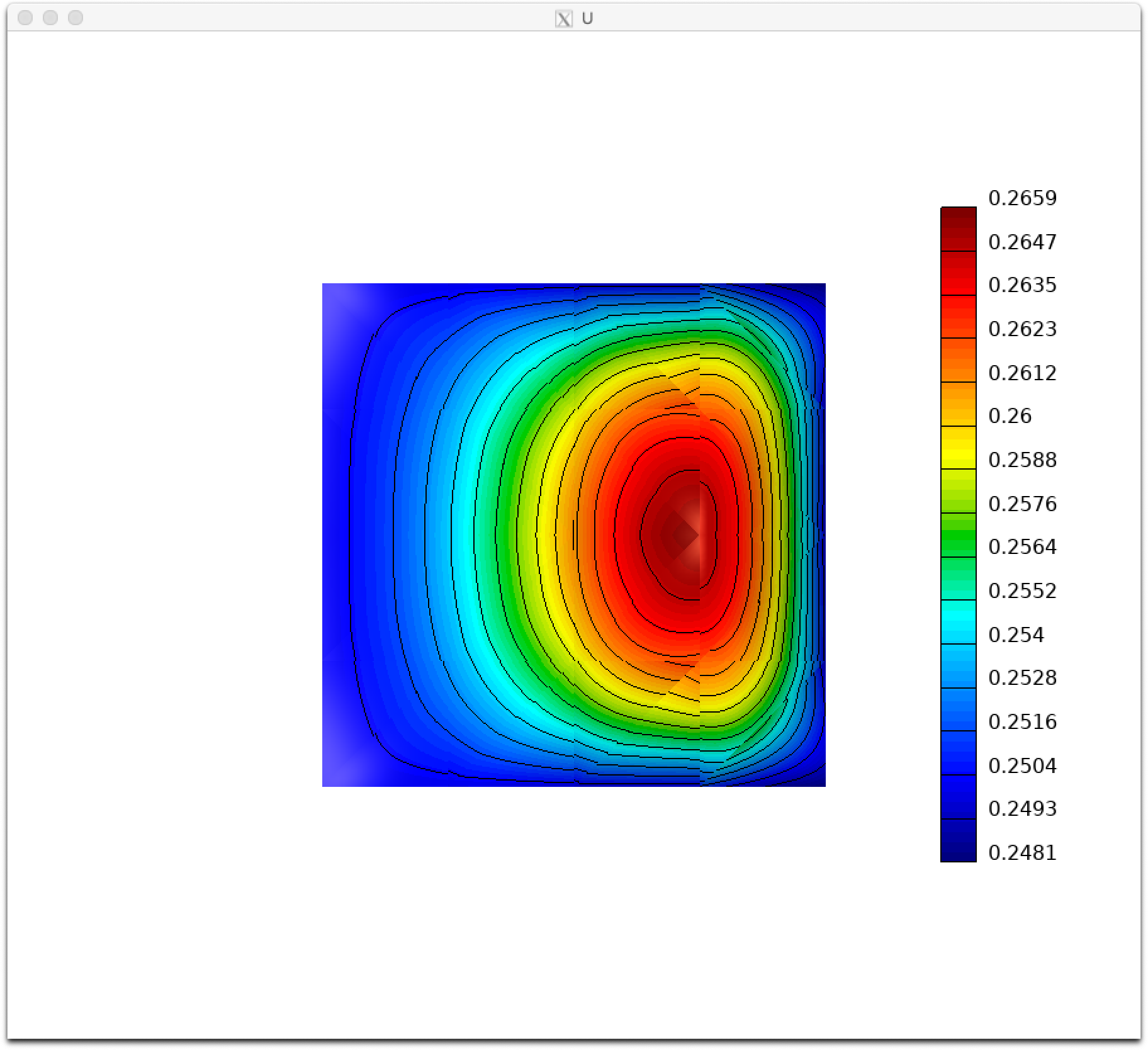}{Initial mesh}{}{$\psi$}{0.249}{0.266};
 \drawContour{xshift=5.25cm,yshift=6.75cm}{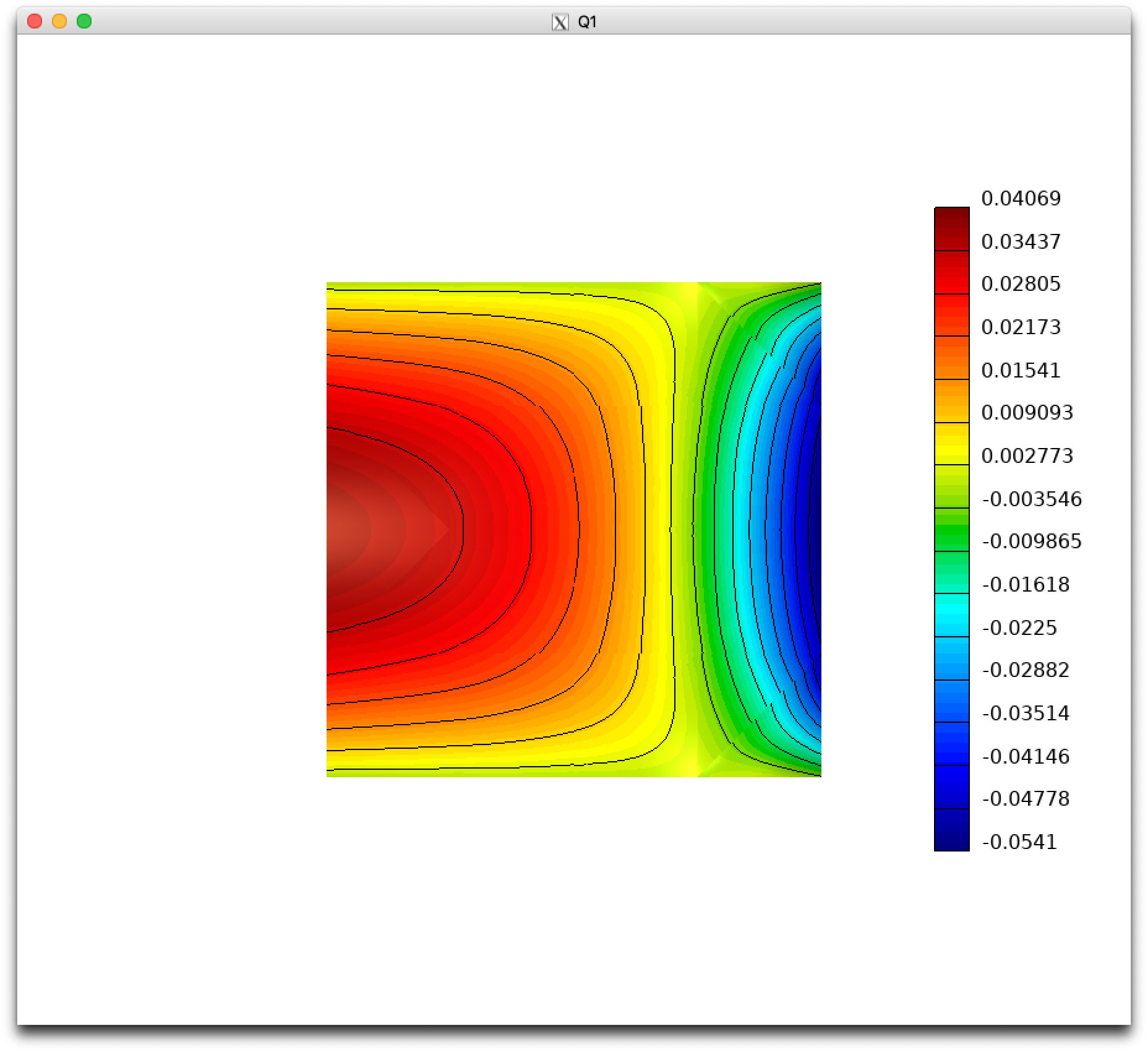}{Initial mesh}{}{$\psi_r/r$}{-0.0478}{0.0407};
 \drawContour{xshift=10.5cm,yshift=6.75cm}{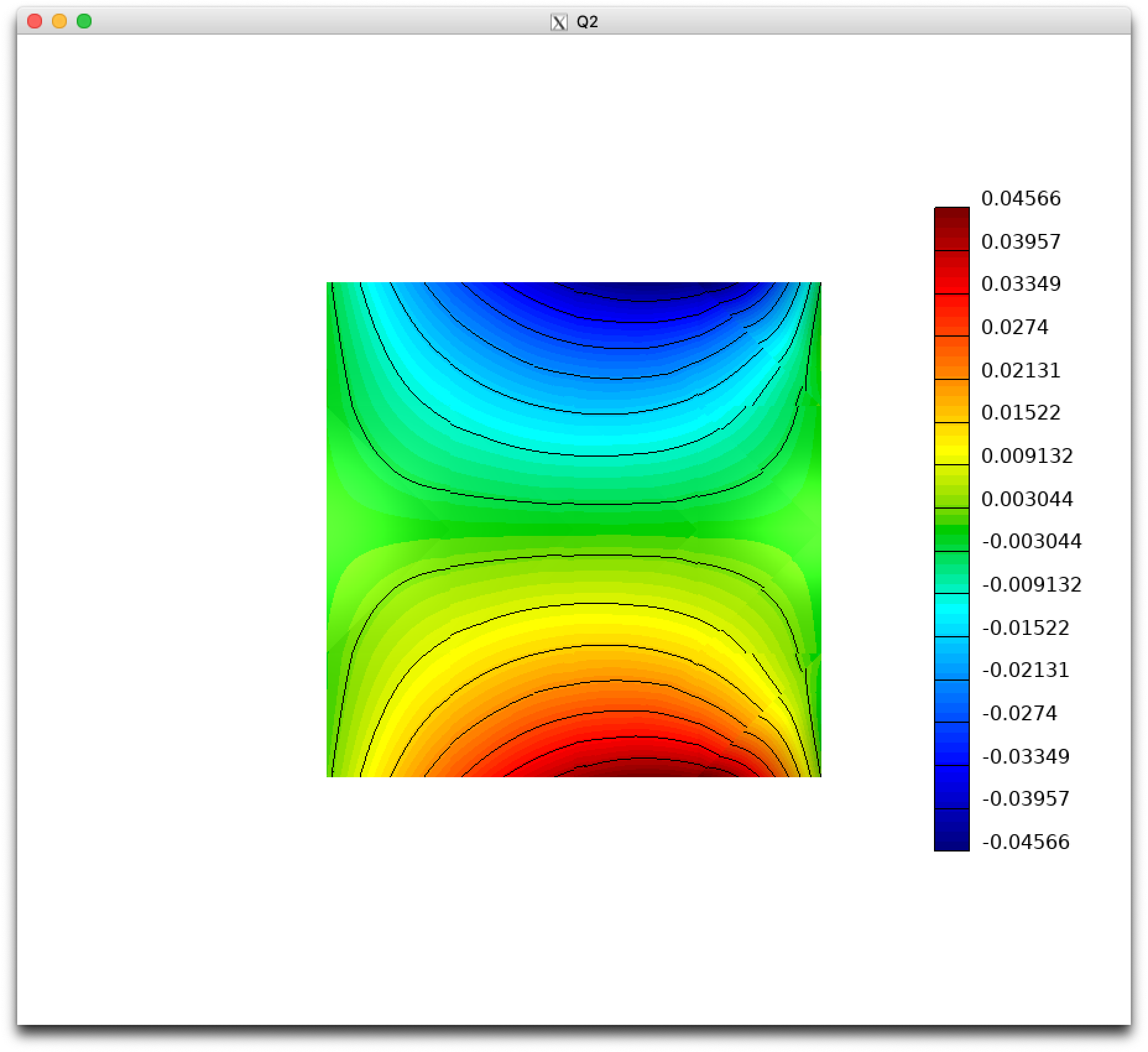}{Initial mesh}{}{$\psi_z/r$}{-0.0396}{0.0457};
 \drawContour{xshift= 0.0cm,yshift=0.5cm}{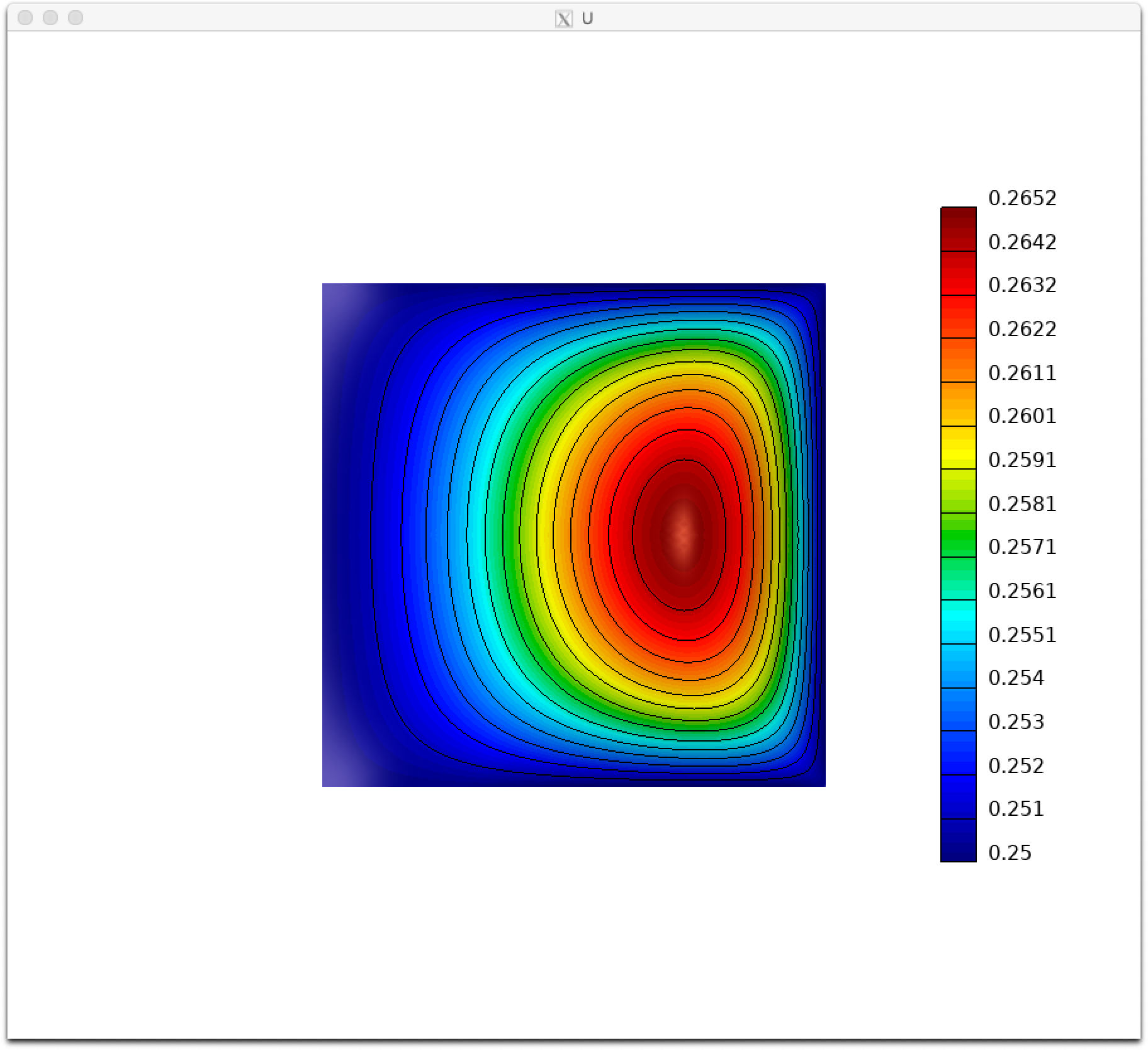}{After 4 refinement}{}{$\psi$}{0.251}{0.265};
 \drawContour{xshift=5.25cm,yshift=0.5cm}{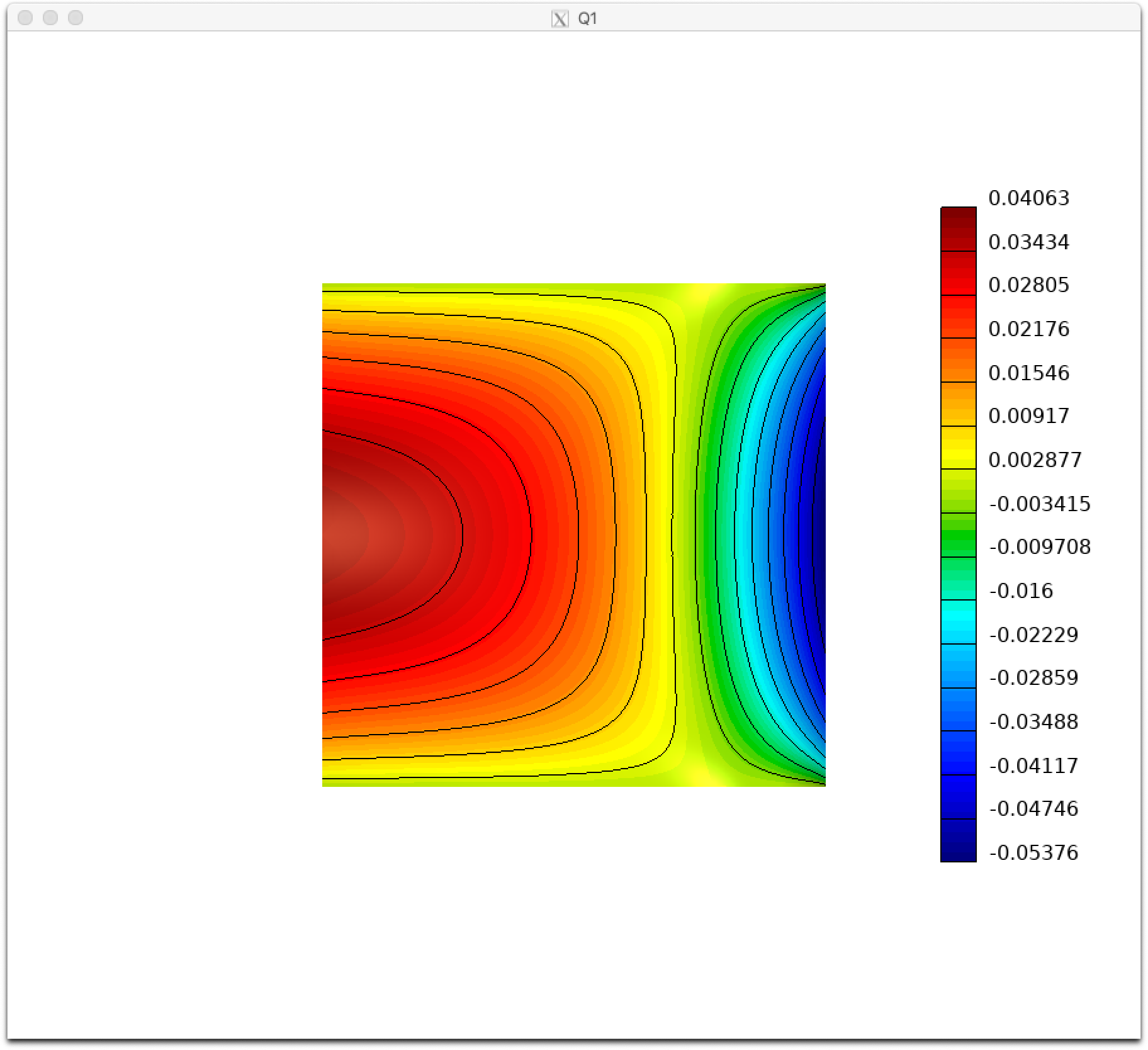}{After 4 refinement}{}{$\psi_r/r$}{-0.0475}{0.0406};
 \drawContour{xshift=10.5cm,yshift=0.5cm}{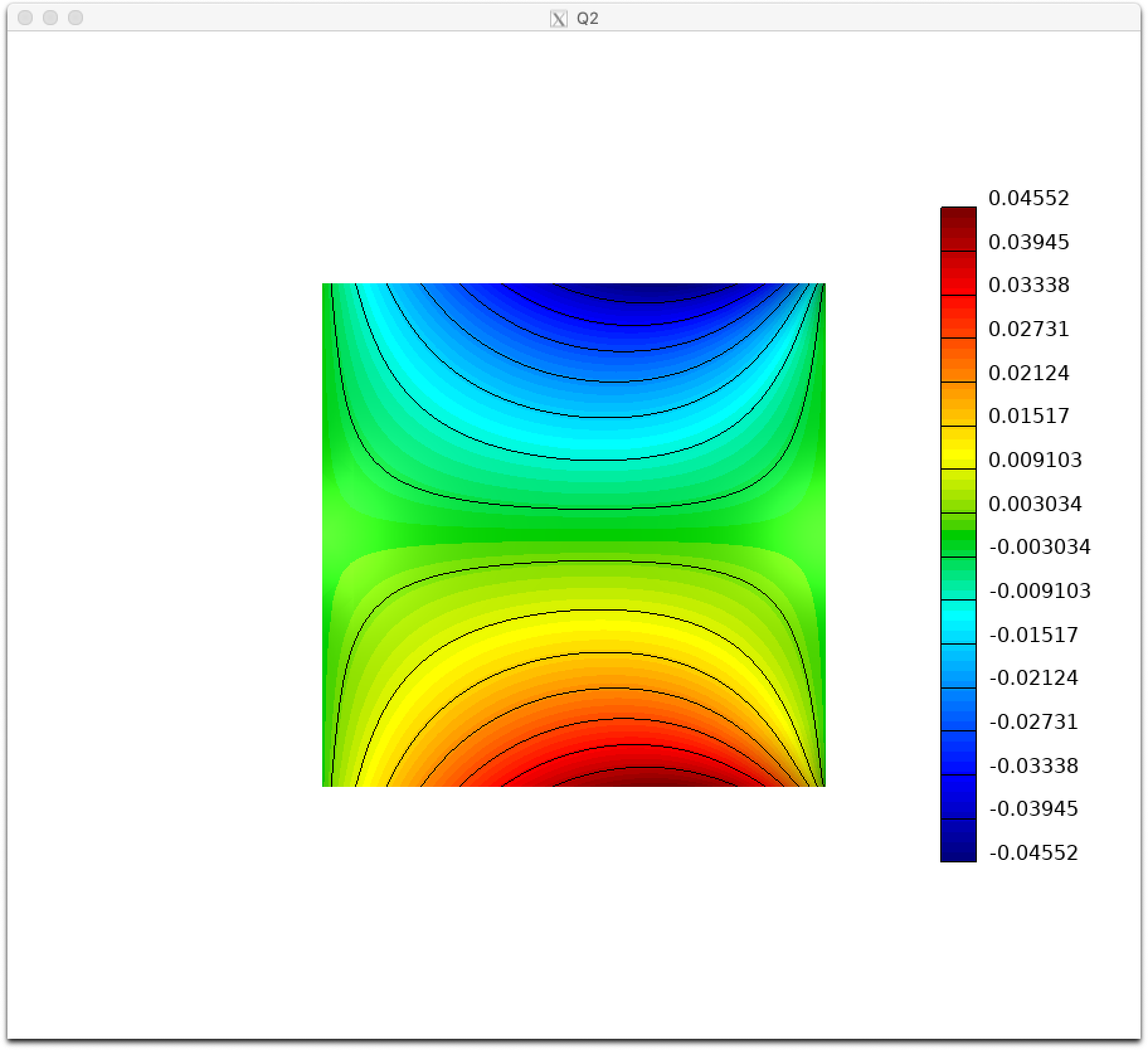}{After 4 refinement}{}{$\psi_z/r$}{-0.0395}{0.0455};
\end{tikzpicture}
}
\end{center}
\caption{Numerical solutions and their derivatives for the rectangular geometry case in Section \ref{sec:amr_test1}.
The first row shows the results on the initial mesh.
The second row shows the results after four AMR iterations.
}
\label{fig:amr1}
\end{figure}
}


\subsubsection{D-shaped geometry}\label{sec:amr_dshape}
To further confirm the effectiveness of AMR on problems involving complex geometry, we consider the same nonlinear test in D-shaped geometry as described in Section \ref{sec:nonlinear_test}. For this test, we set $\textrm{atol}_{\textrm{amr} }=$\pzc{1e-6} and use linear polynomial space $U_h^1$ as the discrete trial space for the DPG scheme. 

The initial mesh and the convergence history of AMR and uniform refinement are presented in Figure \ref{fig:amr_dshape}. We observe that, compared with the uniform refinement, the numerical error on AMR is much smaller and decays faster. This example implies the efficiency and effectiveness of our AMR strategy on problems with complex geometry.

{
\newcommand{\figWidth}{5.0cm}
\newcommand{\figWidthb}{4.9cm}
\newcommand{\trimfig}[2]{\trimh{#1}{#2}{.25}{.2}{.2}{.2}}
\newcommand{\figWidtha}{3.6cm}
\newcommand{\trimfiga}[2]{\trimh{#1}{#2}{.0}{.0}{.0}{.0}}
\begin{figure}[htb]
\begin{center}
\begin{tikzpicture}[scale=1]
  \useasboundingbox (0,.5) rectangle (12.,5.);  
  \draw(0.5,-.3) node[anchor=south west,xshift=-4pt,yshift=+0pt] {\trimfig{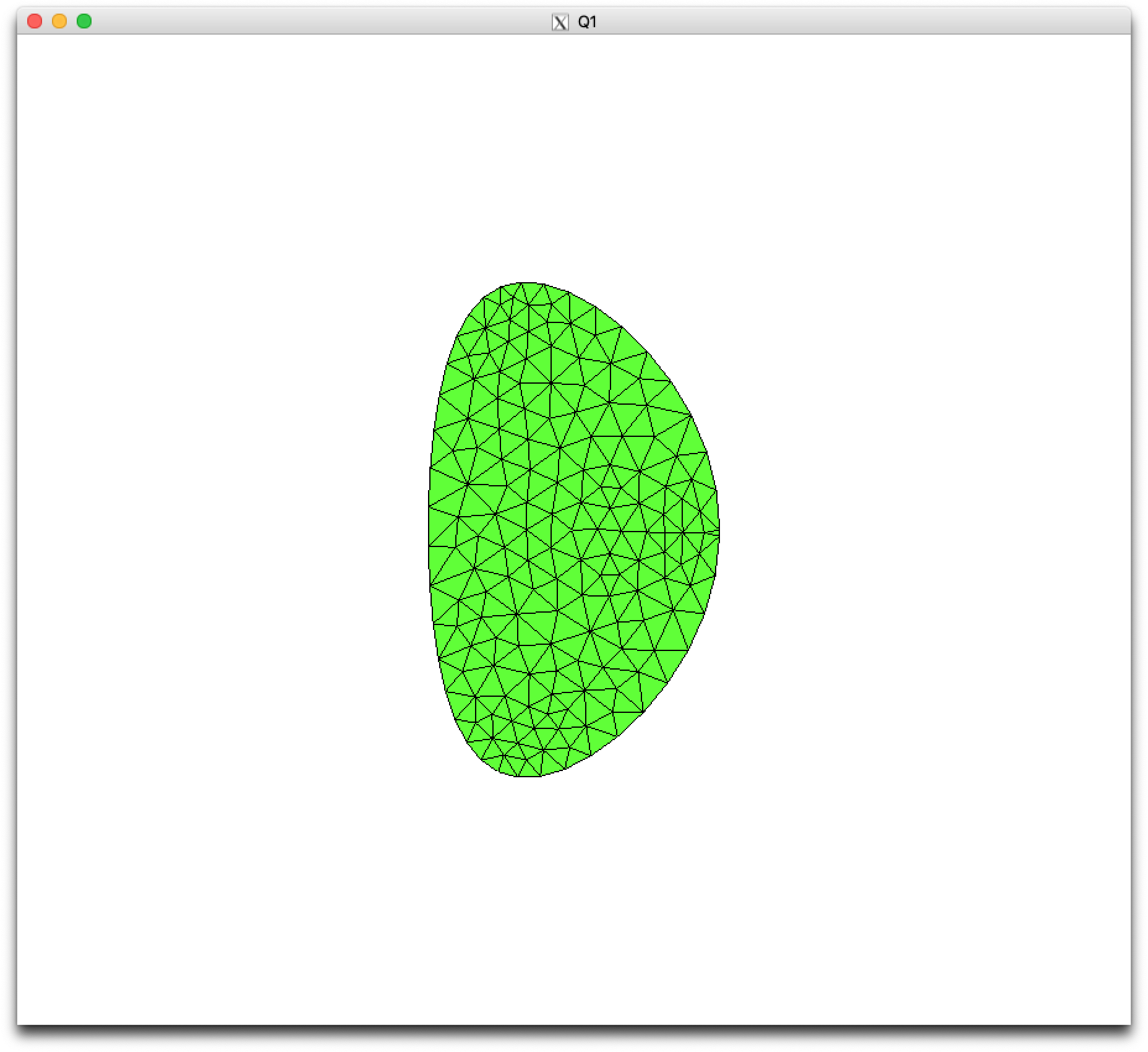}{\figWidth}};
  \draw(6, .3) node[anchor=south west,xshift=-4pt,yshift=+0pt] {\trimfiga{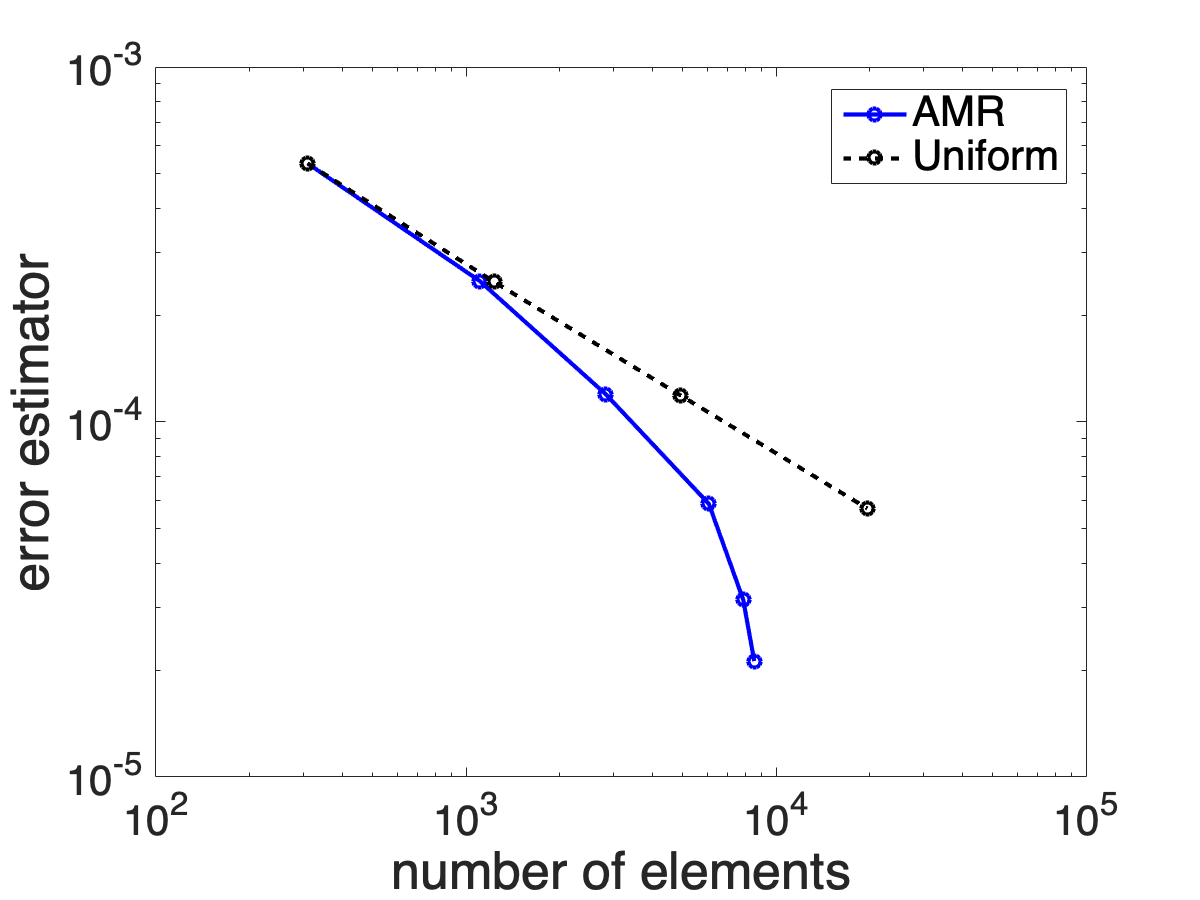}{\figWidtha}};
  \draw(3,4.6) node[draw,fill=white] {\scriptsize Initial mesh};
   \draw(9,4.1) node {\scriptsize Errors vs element number};
%
\end{tikzpicture}
\end{center}
\caption{Left: initial mesh for the test in Section \ref{sec:amr_dshape}.
Right: errors verse total number of elements for both uniform and adaptive refinements. }
\label{fig:amr_dshape}
\end{figure}
}